\documentclass[a4paper,12pt]{article}

\title{Some numerical results in complex differential geometry}
\author{S. K. Donaldson}
\begin{document}

\maketitle
\newcommand{\bC}{\mbox{${\bf C}$}}
\newcommand{\bZ}{\mbox{${\bf Z}$}}
\newcommand{\hdelx}{\mbox{$\hat{\delta_{x}}$}}
\newcommand{\hdely}{\mbox{$\hat{\delta_{y}}$}}
\newcommand{\bR}{\mbox{${\bf R}$}}
\newcommand{\bP}{\mbox{${\bf P}$}}
\newcommand{\dbd}{\mbox{$\overline{\partial}\partial$}}
\newtheorem{prop}{Proposition}
\newtheorem{lem}{Lemma}

\ \ \ \ \ \ \ \ \ \ \ \ \ \ \ {\it Dedicated to Professor Friedrich Hirzebruch}

\section{Introduction}

Far-reaching existence theorems, exemplified by Yau's solution of the Calabi
conjecture for Kahler-Einstein metrics, are a distinctive feature 
of Kahler geometry. There are many open problems---the extension
of Yau's results to the \lq\lq positive'' case and to 
Kahler Ricci solitons;
existence theory for  constant scalar curvature
and extremal metrics; the study of the associated parabolic Ricci flow
and Calabi flow equations---making up  a very active and
challenging research area today. Alongside this, little attention has been
paid to the search for explicit {\it numerical} solutions: that is, good
approximations to the metrics treated by the theory. The only steps
 in this direction seems to be the pioneering and recent work of
 Headrick and Wiseman \cite{kn:HW}. In this article we discuss
  another approach to this question, illustrated by
 some numerical results  for a particular K3 surface $S$; the double cover of
 the projective plane branched over the sextic curve $x^{6}+y^{6}+z^{6}=0$.
The Kahler-Einstein metric $\omega$ on $S$ is characterised by the fact that, after
suitable normalisation, the norm $\vert \theta\vert_{\omega}$ of the standard
holomorphic $2$-form $\theta$  is equal to $1$ at each point. We
find a metric $\omega_{9}'$ on $S$, given explicitly in terms of 26 real
parameters, such that the norm $\vert \theta\vert_{{\omega}_{9}'}$ differs
from $1$ by at most about $1.5\%$ and on average, over $S$, by about
$.11\%$. Thus there are good grounds for thinking that ${\omega}_{9}'$ is
a very reasonable approximation to the Kahler-Einstein metric, which could
be applied to investigate any specific geometric question. The author has
not yet attempted to investigate such applications in this case, except for some discussion of the spectrum of the Laplacian which we
which give in Section 4.   We will focus on
this one example, although it should become clear that the methods could be
applied, in a practical fashion, to other cases. The author has some
other results for  toric surfaces which he hopes to describe
elsewhere. 
 
In Section 1 below we describe a general approach to these approximation
questions for projective algebraic manifolds. This approach is quite different
from that of Headrick and Wiseman. The distinctive features are
\begin{itemize}
\item The use of metrics furnished  by projective embeddings of the manifold.
\item The construction of approximations to the differential-geometric
solutions as limits of iterates of maps
defined by  integration over the manifold.
\end{itemize}
 (It seems likely that these maps, which are dynamical systems on finite-dimensional
 spaces, can be viewed as discrete approximations to the Ricci and Calabi
 flows.)

To expand on the the first item above, consider an ample line bundle $L\rightarrow
X$ over a compact complex manifold $X$, so for large enough integers $k$ the
sections of $L^{k}$ give a projective embedding of $X$. Suppose $G$ is a positive definite
 Hermitian form on the vector space $H^{0}(L^{k})$ of sections. Then there
 is a metric $h$ on the line bundle $L^{k}$ characterised by the fact that
 if $(s_{\alpha})$ is any orthonormal basis of $H^{0}(L^{k})$ then the
 function
    $$ \sum_{\alpha} \vert s_{\alpha}\vert_{h}^{2}$$
    is constant on $X$.
    The curvature of the unitary connection on
    $L^{k}$ associated to $h$ has the form $-i k \omega_{G}$ where
    $\omega_{G}$ is a Kahler metric on $X$ in the class $ 2\pi  c_{1}(L)$.
    This metric can also be viewed as the restriction of the standard
    Fubini-Study metric on the projective space $\bP(H^{0}(L^{k})^{*})$ 
    (defined by the Hermitian form $G$) to the image of $X$ under the 
    projective embedding. Thus we have a way to generate Kahler metrics
    from simple algebraic data---Hermitian forms on finite dimensional
    vector spaces. 
    
    The potential utility of this
     point of view has been advocated over many years by Yau, and a fundamental
     result of Tian \cite{kn:T} shows that this scheme does yield a way to approximate
     any metric. Tian proved that for any metric $\omega$ in the class $2\pi
     c_{1}(L)$ there is a sequence $G_{k}$ of Hermitian forms such that
     $$ \Vert \omega -\omega_{G_{k}}\Vert =O(k^{-2}) , $$
     for any $C^{r}$ norm on $X$. Now the dimension of $H^{0}(L^{k})$ grows
     like $k^{n}$ where $n={\rm dim} X$, so---from a practical point of
     view---one might think that good approximations by these \lq\lq algebraic''
     metrics would require the use of unduly large vector spaces (with
     a Hermitian metric depending on $O(k^{2n})$ real parameters). This
     point can be addressed by a refinement of Tian's result. There is
     a sequence $\tilde{G}_{k}$ such that 
     $$ \Vert \omega-\omega_{\tilde{G}_{k}}\Vert = O(k^{-\nu})$$
     for {\it any } $\nu$. This is explained in the Appendix below.
     In other words, any metric can be very rapidly approximated by algebraic
     ones. (These approximations are in some ways analogous
     to the approximations of a smooth function defined by truncating the
     Fourier series.)  Now suppose we are in a situation where we know (or
     hope)  that a special metric $\omega$ (Kahler-Einstein, constant scalar
     curvature, extremal, Kahler-Ricci solution) exists. How can
     we generate a  sequence of approximations to $\omega$  via Hermitian
     forms $\tilde{G}_{k}$,
     and can we get useful  approximations with values of $k$ which are small
     enough to be manageable in practice? These are the questions we take
     up below.  

          The focus of this article is on the explicit numerical results.
          While there is already in place a considerable quantity of
           rigorous theory to back up
          these methods (and the authors interest in these questions grew
          out of work on abstract existence questions), there are also many
          points where more theory needs to be filled in; but we will not
          dwell on these here. It is also worth mentioning that many of
          these constructions are closely related to ideas in Geometric
          Quantisation Theory, and in particular the asymptotics of the
          classical limit, but we will not say more about this.
         
         Lacking any background in numerical methods, the author is not
         really qualified to comment on the comparison of the results here
         with those of Headrick and Wiseman. In one direction, the end-products
         seem comparable: approximations to Kahler-Einstein metrics on
         very special K3 surfaces with large symmetry groups using standard
         PC's. The author has not been able to directly compare the accuracy
         of the approximations in the two cases, but suspects that the
         best approximation we present here is comparable to the middle
         range of those achieved by Headrick and Wiseman and that by increasing
         $k$ a bit we would achieve something comparable to their best
         approximation. One fundamental disadvantage of our approach is
         that it is limited  to algebraic varieties, whereas Headrick and
         Wiseman are able to vary the Kahler class continuously. One clear
         advantage, at least in  the case studied here, is that we only
         need a few real parameters to store the metrics, as opposed to
         thousands or millions of parameters needed to record the values
         for a lattice approximation.
         
  The author is grateful to Andr\'es Donaldson for instruction in 
  programming. He would also like to mention the crucial part played by
   Appendix A (\lq\lq Elementary programming in QBasic'') of the
  text \cite{kn:A} in getting this project off the ground.

  \section{General theory}
  \subsection{Constant scalar curvature}

We return to the discussion of an ample line bundle $L\rightarrow X$,
 as in the previous section. Thus we consider the relations between
 two different kinds of metric data:
   \begin{itemize}  
     \item Hermitian metrics $G$ on the finite-dimensional complex vector
     space $H^{0}(X;L^{k})$; 
     \item Hermitian metrics $h$ on the line bundle $L^{k}$ such that the compatible
      unitary connection has curvature $-2\pi i k \omega_{h}$, where 
      $\omega_{h}$ is a Kahler
     form on $X$.\end{itemize}
     As a matter of notation, we write $G$ for the metric on $H^{0}(L^{k})$
     and $G^{-1}$ for the induced metric on the dual space. We also use
     the notation $G=(G_{\alpha \beta}), G^{-1}= G^{\alpha \beta}$ when
     working with a basis $(s_{\alpha})$ of $H^{0}(L^{k})$.

     We have two fundamental constructions. Given a metric $h$ on
     $L^{k}$ we let ${\rm Hilb}(h)$ be the Hermitian metric on $H^{0}(L^{k})$;
     $$  \Vert s \Vert_{Hilb(h)}^{2}=  R \int_{X} \vert s\vert^{2}_{h} d\mu_{h} ,
     $$ where $d\mu_{h}$ is the standard volume form $\omega_{h}^{n}/n!$
     and $R$ is the ratio
     $$  R= \frac{{\rm dim} H^{0}(L^{k})}{ {\rm Vol}(X,d\mu_{h})}$$
     (which does not depend on the choice of $h$).
     In the other direction, given a Hermitian metric $G$ on $H^{0}(L^{k})$
     there is a metric $FS(G)$ on $L^{k}$ characterised by the fact that, for
     any orthonormal basis $s_{\alpha}$ of $H^{0}(L^{k})$, we have
     $$      \sum_{\alpha} \vert s_{\alpha}\vert_{FS(H)}^{2}=1, $$
     pointwise on $X$. The Kahler form $\omega_{FS(G)}$ is the restriction
     of the standard Fubini-Study metric on projective space (defined by $G$) to 
     $X$---regarded as a projective
     variety.
Now we say that a pair $(G,h)$ is \lq\lq balanced'' if $G={\rm Hilb}(h)$ and
$h=FS(G)$. (This terminology was introduced in \cite{kn:D1}: the notion had
been considered before by Zhang\cite{kn:Zh} and Luo \cite{kn:Luo}.) Equally, since
in this situation either of $h$ and $H$ determine the other, we can speak of a
metric $h$ on $L^{k}$ or a metric $G$ on $H^{0}(L^{k})$ being \lq\lq balanced''.
Let $M$ denote the set of hermitian metrics on $H^{0}(L^{k})$ and define
$T:M\rightarrow M$ by
$$  T(G)= {\rm Hilb}(FS(G)). $$
Thus, by definition, a balanced metric is a fixed point of $T$.
Now we have, from \cite{kn:D2},
\begin{prop}
Suppose that the automorphism group of the pair $(X,L)$ is discrete. 
If a balanced metric $G$ in $M$ exists then, for any point
$G_{0}\in M$, the sequence $T^{r}(G_{0})$ converges to $G$ as $r$ tends
to infinity.

\end{prop}

(In the original, preprint, version 
 of \cite{kn:D2} this
convergence was raised as a question. In the version submitted for publication
 the convergence property was stated as a fact, but without
a detailed proof. Meanwhile, independently, Y. Sano \cite{kn:S} supplied the detailed
argument, in reponse to the original version of \cite{kn:D2}.) 

We can spell out more explicitly the definition of the map $T$. Changing
point of view slightly,  let
$z_{\alpha}$ be standard homogeneous co-ordinates on $\bC\bP^{N}$ and 
$X\subset \bC\bP^{N}$ be a projective variety. Start with a positive definite
Hermitian matrix  $G_{\alpha \beta}$ and form the inverse $G^{\alpha \beta}$.
For $z\in \bC^{N+1}$ set
$$ D(z)= \sum G^{\alpha \beta} z_{\alpha} \overline{z}_{\beta}. $$
 Then the quotients
$$  f_{\gamma \delta} = \frac{z_{\gamma} \overline{z}_{\delta}}{D(z)}, $$
are homogeneous of degree $0$ and can be regarded as functions on
$\bC\bP^{N}$, and so on $X$. Then the map $T$ is
$$  (T(G))_{\gamma \delta} = R \int_{X} f_{\gamma \delta} d\mu, $$
where $d\mu$ is the volume form induced by the Fubini-Study metric of
$G$. Notice that the constant $R$ is chosen so that $G^{\gamma \delta} T(G)_{\gamma
\delta}= N+1$. The geometry is unaffected by rescaling the metric $G$ so
in practise we work with metrics normalised up to scale in some convenient but arbitrary
way.  

So far the parameter $k$ has been fixed; we now consider the effect of
increasing $k$. 
The main result of \cite{kn:D1} relates balanced metrics
  to constant scalar curvature metrics on $X$.
\begin{prop}
   Suppose that the automorphism group of the pair $(X,L)$ is discrete.
   If $X$ has a constant scalar curvature Kahler metric $\omega$ in the
   class $ 2\pi c_{1}(L)$ then for large enough $k$ there is a unique balanced
   metric on $L^{k}$ inducing a Kahler metric $ k \omega_{k}$ on $X$, and
   $\omega_{k}\rightarrow \omega$ as $k\rightarrow \infty$. Conversely,
   if there are balanced metrics on $L^{k}$ for all large $k$ and the sequence
   $\omega_{k}$
   converges then the limit has constant scalar curvature.
\end{prop}

Taken together, these two results give a procedure for finding numerical
approximations to constant
scalar curvature metrics. We choose a sufficiently large value of $k$ and
then compute the iterate $T^{r}(G_{0})$ for some convenient initial metric
$G_{0}$ on $H^{0}(L^{k})$. If $k$ is sufficiently large then the limit
as $r$ tends to infinity is a good approximation to the differential geometric
solution.

A crucial ingredient in the proof of Proposition 2 is the Tian-Yau-Zelditch-Lu expansion
for the \lq\lq density of states'' function. If $s_{\alpha}$ is an orthonormal
basis for $H^{0}(L^{k})$ with respect to the standard $L^{2}$ metric we
set $\rho_{k}= \sum \vert s_{\alpha}\vert^{2}$; a function on $X$ which
does not depend on the choice of basis. Then if we take a {\it fixed} metric
on $L$ and form the sequence of functions $\rho_{k}$ with the induced metrics
on $L^{k}$ we have
$$\rho_{k} \sim k^{n} + a_{1} k^{n-1} + a_{2} k^{n-2} + \dots, $$
where the $a_{i}$ are local invariants and $a_{1}$ is $1/2\pi$ times the scalar
curvature. 

\subsubsection{A toy example}
We take $X$ to be the Riemann sphere $\bC\bP^{1}$, $L$ to be the line
bundle ${\cal O}(1)$. (Strictly this example does not fit into
the framework above, since the pair has a continuous automorphism group
$SL(2,\bC)$. However, the theory can undoubtedly be extended to relax the
condition on the automorphisms, in the manner of \cite{kn:Mab}, so we will ignore
this technicality here.) A basis of $H^{0}(L^{k})$ is given by
$1,x,x^2,\dots x^k$ where $x$ is a standard co-ordinate on $\bC$. We restrict attention to $S^{1}$-invariant metrics,
for the  $S^{1}$ action $x\mapsto e^{i\theta}x $ on the sphere. The invariant
metrics are represented by diagonal matrices in our basis, specified by
the $k+1$ diagonal entries, thus
$$D= \sum_{p=0}^{k} a_{p}\vert x\vert^{2p}. $$
(Notice that the $a_{p}$ are really the entries of the metric $G^{-1}$ on
the dual space: in practice it is easier to work with this rather than
$G$.)
The round metric on the sphere is given by $a_{p}=\left( \begin{array}{c}
 k\\ p \end{array}\right)$, when $
D= (1+\vert x\vert^2)^k$. Due to the symmetry we know that for any $k$
the balanced metric is a standard round metric on the sphere. 

To make things even simpler, we can consider metrics invariant under
the inversion $x\mapsto x^{-1}$, so $a_{p}= a_{k-p}$. Since the geometry
is unaffected by a overall
scaling $a_{p}\mapsto Ca_{p}$ we have just $\lfloor k/2\rfloor+1$ essential real
parameters. In the first case, when $k=2$, one can evaluate the integrals
in the definition of $T$ explicitly using elementary calculus. If we normalise
so that $a_{0}=a_{2}=1/2$ and write $a_{1}=s$ the map is represented by
the function
$$  \tau(s)=  \frac{ s \cosh^{-1}(s) + \sqrt{s^{2}-1} (s^{2}-2)}{ 2s\sqrt{s^{2}-1}
- 2 \cosh^{-1}(s)} . $$
The reader who plots this function will immediately see that the iterates
$\tau^{r}(s_{0})$ do indeed converge rapidly to the fixed point $s=1$.
(The formula as written is valid in the range $s>1$ but has an obvious
continuation to $s\leq 1$.)

We now consider the case $k=6$, and evaluate the integrals numerically.
A typical sequence of iterates is indicated in the next table.

$$\begin{tabular}{|l||l|l|l|l|}\hline
r& $a_{0}$&$a_{1}$&$a_{2}$&$a_{3}$\\ \hline
0&.018&.495&4.5&54\\ \hline
1&.02833&.8539&11.04&40.16\\ \hline
2&.03923&1.268&13.38&34.62\\ \hline
3&.05331&1.645&14.39&31.81\\  \hline
4&.07150&1.987&14.90&30.09\\  \hline
10&.2384&3.493&15.57&25.40\\ \hline
20&.7365&5.400&15.26&21.20\\ \hline
30&.9488&5.895&15.05&20.21\\ \hline
40&.991&5.983&15.01&20.03\\ \hline
$\infty$&1&6&15&20\\ \hline   
\end{tabular}$$

This confirms the convergence that the theory predicts, and similar resuts
are obtained whatever initial values for $a_{0}, a_{1}, a_{2}, a_{3}$
are used. The choice of the particular initial values here is made because
in this case all the metrics in the sequence can be represented as surfaces
of revolution in $\bR^3$, so one obtains a vivid representation of the
evolution of the geometry of the surface through the sequence. The initial
values give, roughly speaking, a connected sum of spheres joined by a  small
neck. The first application of $T$ stretches the surface into a long \lq\lq
sausage'', which becomes convex after four more application of $T$, and thereafter the sausage slowly shrinks in length to approach
a round sphere.

In this example we see that the convergence to the limit is, while steady,
 quite slow. Let $\epsilon_{i}(r)$ be the difference between $a_{i}(r)$
 and the limiting value $a_{i}(\infty)$. From standard general theory we know that there is a constant
 $\sigma\in (0,1)$ associated to the problem such that for almost all initial
 conditions 
 $$ \epsilon_{i}(r) \sim c_{i} \sigma^{r}
$$ as $r\rightarrow \infty$. This constant $\sigma$ is just the largest
eigenvalue of the derivative of $T$ at the fixed point, and the vector
$c_{i}$ is an associated eigenvector. Analysing the sequence,
 one can see that in this case $\sigma$ is about .8. We also see numerically that
 the corresponding eigenvector has entries approximately 
 $$   (1,-2,-1,-4). $$
 The fact that one gets this eigenvector is easily  explained by the
 $SO(3)$ invariance of the problem.  The eigenvector
  corresponds to the second spherical
 harmonic on the sphere. This is given by
      $$ f_{2}= \left(\frac{1-\vert x\vert^2}{1+\vert x\vert^2}\right)^2.
      $$
      Thus $$f_{2}= \frac{(1-\vert x\vert^2)^{2}(1+\vert x\vert^2)^4}{(1+\vert
 x\vert^2)^6},$$ 
 and $$ (1-\vert x\vert^2)^{2}(1+\vert x\vert^{2})^{4}= 1-2\vert x\vert^{2}-
 \vert x\vert^{4} - 4\vert x\vert^{6} - \vert x \vert^{8} -2 \vert x\vert^{10}+
 \vert x\vert^{12}. $$
 
 A general point to make here is that while the convergence of the sequence
 $T^{r}$ is slow (since .8 is not much less than  1) it is easy, using this standard
 analysis of the linearisation, to define  much more rapidly convergent
 sequences, and of course the same remarks hold in more complicated examples.

This toy model is the only example of the $T$-iteration which we will consider
in this paper. The picture extends readily to the case of {\it toric
 varieties} in higher dimensions  which, as mentioned before, the author
 hopes to take up elsewhere.

\subsubsection{Extremal metrics}

  As we stated above, the results of \cite{kn:D1} and \cite{kn:D2} are limited, strictly, to
  the case where the pair $(X,L)$ does not have continuous automorphisms,
  but we anticipate that the theory can be extended to remove this condition.
  In this regard, we point out here that there is a straightforward modification
  which can be expected to produce numerical approximations to {\it extremal
  metrics} in the sense of Calabi \cite{kn:C}. Suppose that $(X,L)$ has a conected
  Lie
  group of automorphisms and fix a maximal compact subgroup $K$.
  Then
    $K$ acts on $H^{0}(X,L^k)$ and we restrict to $K$-invariant
    metrics. Then it may happen that the sequence $T^{r}(G_{0})$ does not
    converge but that there is a sequence $g_{r}$ in the {\it complexification}
    of $K$ such that $g_{r} ( T^{r}(G_{0}))$ does converge,  as
 $r\rightarrow \infty$. In this case, taking the limits, we again get a
  sequence of preferred metrics $\omega_{k}$ which we expect to converge
  to an extremal metric on $X$ as $k\rightarrow \infty$. While the theory
  here needs to be filled in, the procedure works effectively in the 
  examples
  that the author has studied numerically.  

\subsection{Calabi-Yau metrics}

We now turn to the main topic of this paper: the case when the metric one
wants to approximate is Ricci-flat. While this case could be treated within
the constant scalar curvature theory outlined above, there is a slightly
different, and simpler approach. Let $X,L$ be as before and suppose given
a fixed smooth volume form $d\nu$ on $X$. A fundamental result of Yau asserts
that there is a unique Kahler metric $\omega$ on $X$ in the class 
$ 2\pi c_{1}(L)$ which realises $d\nu$ as its volume form $\omega^n/n!$. 
We propose a method for finding numerical approximations to this metric.
The most interesting case is when $X$ is a Calabi-Yau manifold, so has
a nowhere-vanishing holomorphic $n$-form $\theta$. Then the metric with
volume form $ d\nu= (-i)^{n}  \theta \wedge \overline{\theta}$ is Ricci-flat.

In a nutshell, we modify the definitions of the previous section by using
the fixed volume form $d\nu$ in place of the Fubini-Study volume $d\mu$.
Thus we say that a Hermitian metric ${G}$ on $H^{0}(X;L^{k})$ is 
$\nu$-{\it balanced} if $G={\rm Hilb}_{\nu}(FS(G))$ where
 $$ \Vert s\Vert^{2}_{\rm Hilb_{\nu}(h)} = R  \int \vert s \vert_{h}^{2} d\nu. $$
The existence of balanced metrics in this context is due to Bourguignon,
Li and Yau \cite{kn:BLY}, see also the recent paper \cite{kn:Ar}. It is closely related to well-known results
relating moment maps to Geometric Invariant theory, and similar extensions
were considered by Milson and Zombro \cite{kn:MZ}.
In fact suppose that $\nu$ is any
positive
Radon measure on $\bC\bP^{n}$. We say that a metric  
$G_{\alpha \beta}$
 on $\bC^{N+1}$
 is $\nu$-{\it-balanced} if
 $$ G_{\gamma \delta}= R \int_{\bC\bP^{N}}
  \frac{z_{\gamma}\overline{z}_{\delta}}{D(z)} d\nu. $$
   
    When $\nu$ is derived from a smooth volume
   form supported on a projective subvariety $X\subset \bC\bP^{N}$ this
   reproduces the previous definition. But we can also consider other measures,
   in particular, sums of point masses (or, in other words, the case when
   subvariety has dimension $0$). Let us suppose that the measure $\nu$
   satisfies one of the following two conditions:
   \begin{enumerate}
   \item For any non-trivial linear function $\lambda$ on $\bC^{N+1}$ the
   function $\log(\frac{\vert \lambda(z)\vert}{\vert z\vert})$ on $\bC\bP^{N}$ is
   $\nu$-integrable.
   \item $\nu$ is a sum of point masses supported on a finite set $Z$ and for any projective subspace
   $P\subset \bC\bP^{N}$ we have
   $$    \frac{\nu(Z\cap P)}{{\rm dim}\ P +1} < \frac{\nu(Z)}{N+1}. $$
   \end{enumerate}
   
   (The definition in the first item uses a metric on $\bC^{N+1}$ but is
   clear that the condition does not depend on this choice.) It is easy
   to see that if $\nu$ is a smooth volume form on a subvariety $X$ which
   does not lie in any proper projective subspace then the first hypothesis
   holds.

   Now we have
   \begin{prop}
   If $\nu$ is a positive Radon measure on $\bC\bP^{N}$ which satisfies
   either condition (1) or (2) above then there is a  $\nu$-balanced
   metric on $\bC^{N+1}$, and this is unique up to scale.
   \end{prop}
   
    We recall the proof briefly. Note that the space of metrics $M$ is
    the symmetric space $GL(N+1, \bC)/U(N)$ and there is a standard notion
    of geodesics in $M$; the images of analytic $1$-parameter subgroups
    in $GL(N+1,\bC)$. For any non-zero vector $z\in \bC^{N+1}$
    we let $\psi_{z}$ be the function
    $$  \psi_{z}(G)= \log \vert z \vert^{2}_{G^{-1}} +\frac{1}{N+1} \log \det G $$
    on $M$. The key point is that $\psi_{z}$ is convex on all geodesics.
     Changing
    $z$ by a scalar multiple only changes $\psi_{z}$ by the addition of
    a constant. Now given our Radon measure $\nu$ we set, with some abuse
    of notation,
    $$ \Psi_{\nu}(G)= \int_{\bC\bP^{N}} \psi_{z} d\nu(z). $$
    This is defined up to the addition of a constant; for example we can
    define $\psi_{z}$ for $z\in\bC\bP^{N}$ by taking the lift to a vector
    in $\bC^{N+1}$ of length $1$ with respect to some chosen reference
    metric.  Now $\Psi_{\nu}$ is also convex, being a positive linear combination
    of convex functions. It is easy to check that a metric $G$ is $\nu$-balanced
    if and only if it is a minimum of $\Psi_{\nu}$. Such a minimum will
    exist so long as $\Psi_{\nu}$ is a {\it proper} function on $M$, and the convexity yields uniqueness. In turn, $\Psi_{\nu}$ is
    proper on $M$ if and only if it is proper on each geodesic. So what we have to verify is that---under either of the hypotheses
    (1), (2)---for each geodesic ray $G_{t}$ we have
    $\Psi_{\nu}(G_{t})\rightarrow \infty$ as $t\rightarrow \infty$. 
     Using the $GL(N+1,\bC)$-invariance of the problem it suffices to consider
     a geodesic of the form
     $$   G_{t} = {\rm diag} ( e^{\lambda_{\alpha} t} ), $$
     where $\sum \lambda_{\alpha}=0$ and $\lambda_{0}\geq\dots \geq \lambda_{N}$. 
     Then $$\Psi_{\nu}(G_{t})= \int_{\bC\bP^{N}} \log\left( \sum e^{\lambda_{\alpha}
     t} \vert x_{\alpha}\vert^{2}\right) d\nu_{x}. $$
     Consider first case (1). The first coefficient $\lambda_{0}$ must
     be positive and we have
     $$ \log(\sum e^{\lambda_{\alpha} t} \vert x_{\alpha}\vert^{2} \geq
     \log (e^{\lambda_{0} t } \vert x_{0} \vert^{2})= \lambda_{0} t + \log
     \vert x_{0}\vert^{2}. $$
     By the hypothesis the term on the right hand side is $\nu$-integrable
     and so
     $$ \Psi_{\nu}(G_{t})\geq \lambda_{0} t \nu(\bC\bP^{N})+ {\rm Const.}$$
     and we see that $\Psi_{\nu}(G_{t})\rightarrow \infty$ as required.
     In case (2) we write $\nu=\sum \nu_{i} \delta_{x^{(i)}}$, for points
     $x^{(i)}$ in $\bC\bP^{N}$. Let
     $$ \alpha(i)= {\rm min} \{ \alpha \vert x^{(i)}_{\alpha}\neq 0 \}.$$
     Then $\Psi_{\nu}(G_{t})\sim c t $ as $t\rightarrow \infty$ where
     $$  c= \sum \lambda_{\alpha(i)} \nu_{i}. $$
     Elementary arguments, essentially the same as in usual Geometric Invariant
     Theory discussion in \cite{kn:New}, \cite{kn:Mum}, show that the condition that
     $c>0$, for all geodesics, is equivalent to the hypothesis (2).

     Thus we know that, under the very mild hypotheses (1) or (2), $\nu$-balanced
     metrics exist. To find them, we make the obviuous modification to
     the algorithm of the previous section. We define $T_{\nu}: M\rightarrow
     M$ in just the same way as $T$ but using the measure $\nu$.
     That is, starting with a matrix $G$, we set
     $$ T_{\nu}(G)_{\gamma \delta} = R \int_{\bC\bP^{N}}\frac{z_{\gamma} 
     \overline{z}_{\delta}}{D(z)}
     d\nu(z). $$
     Notice that this is unaffected by rescaling $\nu$. 
     \begin{prop}
     Suppose $\nu$ satisfies either hypothesis (1) or (2). Then for any
     initial metric $G_{0}$ the sequence $T_{\nu}^{r}(G_{0})$ converges to the
     $\nu$-balanced metric as $r\rightarrow \infty$.
     \end{prop}
        
     The proof is similar to that for the map $T$, but more elementary.
     We show that $T_{\nu}$ decreases the function $\Psi_{\nu}$; then the
     conclusion follows from the properness of $\Psi_{\nu}$. To make the
     notation simpler we will treat the case of point masses, the other
     case being essentially the same. To prove the inequality
     $\Psi_{\nu}(T_{\nu}(G))\leq \Psi_{\nu}(G)$ we can without loss of generality
     suppose that $G$ is the metric given by the identity matrix. We can also
     suppose that the total mass $\sum \nu_{i}$ is $1$.
      We choose representatives $z^{(i)}$ in $\bC^{N+1}$ with 
      $\vert z^{(i)}\vert_{G^{-1}}=1$.
     Then $$\Psi_{\nu}(G)= \sum \nu_{i} \log \vert x^{(i)} \vert_{G^{-1}}^{2}
     +\frac{1}{N+1} \log (1) = 0 .$$
     We treat the two terms in the definition of $\Psi_{\nu}$:
     $$ \Psi_{\nu}(T(G))= \sum \nu_{i} \vert z^{(i)}\vert^{2}_{T(G)^{-1}} +
     \frac{1}{N+1} \log \det T(G), $$
      separately.
    By the concavity of the logarithm function
    $$ \sum_{i} \nu_{i} \log \vert z^{(i)}\vert^{2}_{T(G)^{-1}} \leq \log \left(
    \sum \nu_{i} \vert z^{(i)}\vert^{2}_{T(G)^{-1}}\right). $$
    Now $$ \sum_{i}\nu_{i} \vert z\vert^{2}_{T(G)^{-1}} = \sum_{i,\alpha,\beta}
    \nu_{i} z^{(i)}_{\alpha}\overline{z}^{(i)}_{\beta} (T(G))^{\alpha \beta},
    $$ but this is
    $$ \sum_{\alpha \beta} \frac{1}{N+1} T(G)_{\alpha \beta} T(G)^{\alpha
    \beta}$$
    which is $1$, since $T(G)_{\alpha \beta}, T(G)^{\alpha \beta}$ 
     are inverse matrices. So
    the first term in the definition of $\Psi_{\nu}$ is less than or equal
    to $0$. For the second term, the 
     arithmetic-geometric mean inequality for the eigenvalues gives
    $$ \frac{1}{N+1} \log \det T(G) \leq \log \left( \frac{{\rm Tr}(T(G))}{N+1}\right)
    $$
    and the term on the right is zero, since
    $$ {\rm Tr}(T(G))= \sum_{i}\nu_{i} \vert z_{i}\vert_{G^{-1}}^{2} =\sum
    \nu_{i}= 1. $$
         
     Putting the two terms together we have $\Psi_{\nu}(T(G) \leq 0 $, as required.

    \
    
    To sum up, given a volume form $\nu$ on our algebraic variety $X$ we
    have, for each $k$, an algorithm for finding the $\nu$-balanced metric
    on $H^{0}(X;L^{k})$. Moreover, this is robust in the sense that if
    we approximate $\nu$ by another measure $\nu^{*}$ which is a sum of point masses---as we have to do in
    numerical integration---the numerical algorithm defined by $T_{\nu^{*}}$
    will converge to a $\nu^{*}$-balanced metric provided only that $\nu^{*}$
    satisfies the very mild condition (2) (which will happen for any reasonable
    approximation). Taking the restriction of the Fubini-Study metric and
    scaling by $k^{-1}$
    we get a $\nu$-balanced Kahler metric $\omega_{k,\nu}$ on $X$.
     Now let $k$ tend to infinity. We expect that the result
    corresponding to Proposition 2 is true, so that the kahler metrics
    $k^{-1} \omega_{k,\nu}$ converge to the kahler metric with volume form
    $\nu$. While this is, again, a piece of theory that needs to be 
    filled in, we will assume it is so for the rest of this paper. (In
    one direction, it is not hard to see that if the $k^{-1} \omega_{k,\nu}$
    converge the limit must have volume form $\nu$. We expect that the
    harder converse can be proved by arguments similar to those in \cite{kn:D1}.)
    Thus we have another procedure for finding numerical approximations
    to Calabi-Yau metrics. Experimentally at least, this converges  more 
    quickly than the more general \lq\lq $T$-algorithm'', which is probably
   related to  the fact that the constant scalar curvature equation
    is of fourth order in the Kahler potential while the Calabi-Yau equation
    is a second order Monge-Ampere equation.
     
     To give a toy example, consider again the $S^{1}$ invariant metrics
     on $\bC\bP^{1}$ and the sections of ${\cal O}(6)$. Let $\nu$ be the
     volume form of the standard round metric. Then, with the same starting
     point as before, we obtain the  iterates of  $T_{\nu}$ shown in the
     following table.
     
    $$ \begin{tabular}{|l||l|l|l|l|}\hline
     $r$ & $a_{0}$&$ a_{1}$&$a_{2}$&$a_{3}$\\ \hline
    0& .018&.5&4.5&54\\ \hline
    1&.1395&2.599&15.20&28.12\\ \hline
    2&.4200&4.420&15.54&23.23\\ \hline
    3&.6920&5.297&15.30&21.41\\ \hline
    4&.8568&5.697&15.14&20.61\\ \hline
    10&.9992&5.998&15.00&20.00\\ \hline
    13& .9999&6.000&15.00&20.00\\ \hline
     \end{tabular}$$
     
      The convergence is much faster. The parameter $\sigma$ governing
      the asymptotic behaviour is now  about .42 (with the same eigenvector).

\subsubsection{Refined approximations}

In either variant of the theory, the balanced metrics or $\nu$-balanced
metrics cannot usually be expected to give very close approximations to the desired
differential-geometric solutions for practical values of $k$. An analysis
of the convergence, using the Zelditch expansion, would probably show that
the convergence is only $O(k^{-1})$ or $O(k^{-2})$. We now return to the
issue raised in the Introduction of finding {\it rapidly convergent} approximations.
We restrict the discussion here to the Calabi-Yau case. 
Thus we suppose that we have a metric $G_{0}$ on $H^{0}(X;L^{k})$ inducing
a Fubini-Study metric $\omega_{G_{0}}$ on $X$, and that the volume form
$d\mu_{G_{0}}= \frac{1}{n!} \omega_{G_{0}}^{n}$ is reasonably close, but not extremely close
to the given volume form $\nu$ (which we assume to be normalised so that
the total volumes are equal). We set $\eta=d\mu_{G}/\nu$; a function on
$X$ which is close to $1$. We would like to \lq\lq refine'' $G_{0}$ to get a better
approximation. This is essentially
a linear problem. In the standard differential geometric theory we would
consider a nearby Kahler metric of the form $\omega_{G_{0}} +i \dbd \phi$. The linearisation
of the volume form map is one half the Laplacian of $(X,\omega_{G_{0}})$ so we would
take
$$\phi = -2 \Delta^{-1} ( \eta-1)$$
 where $\Delta^{-1}$ is the Green's operator.
Then, provided that $\omega_{G_{0}}$ is a good enough approxiation for the
linearisation to be valid, the metric $\omega_{G_{0}}+i\dbd \phi$ would be
a better approximation to the desired Calabi-Yau metric and, iterating
the procedure, we would generate a sequence which converged to that limit.
Thus the question we address here is how to implement a procedure like
this numerically, staying within the class of the \lq\lq algebraic'' Kahler
metrics. 

Let $G$ be any metric on $H^{0}(L^{k})$ and 
$s_{\alpha}$ be  a basis of sections. Then we have functions $f_{\alpha \beta}=(s_{\alpha}, s_{\beta})$
on $X$, where $(\ ,\ )$ denotes the pointwise inner product on $L^{k}$
induced by $G$.  Define \lq\lq $\eta$-
coefficients''
$$  \eta_{\alpha \beta}= R \int_{X} f_{\alpha
\beta} (\eta-1) d\nu. $$
(Notice here the diagonal terms $\eta_{\alpha \alpha}$ have a particularly simple
interpretation, when $G$ is the balanced metric and the basis is orthonormal.
 In that case
  the integral of the positive function
   $f_{\alpha \alpha}$ is 
$ 1/R $ and the numbers $\eta_{\alpha \alpha}$ 
can be viewed as a collection of  mean values of $\eta-1$,
weighted by $ R f_{\alpha \alpha}$.) A natural criterion to define an 
\lq\lq optimal'' algebraic approximation to the differential-geometric
solution is that all the coefficients $\eta_{\alpha \beta}$ vanish.
Thus we say that a metric $G$ near to the balanced metric $G_{0}$ is a
\lq\lq refined approximation'' if this occurs. One could hope to prove
that, for large enough $k$, these refined approximations exist and (with
a suitable interpretation of \lq\lq close'') are unique. Further, it is
reasonable to expect that the resulting sequence of refined approximations
would be very rapidly convergent, in the manner discussed in the Introduction.

 In more invariant terms, we define a vector
space ${\cal H}(L^{k})$ to be the Hermitian forms  on $H^{0}(L^{k})$.
 Then, given $G_{0}$,  we have  maps
$$ \iota: {\cal H}(L^{k})\rightarrow C^{\infty}(X)\ \ \  \pi: C^{\infty}(X)\rightarrow
{\cal H}(L^{k}) $$
defined in terms of a basis by
$$ \iota((a_{\alpha \beta}) = \sum a_{\alpha \beta} f_{\alpha \beta}, $$
and 
$$ \left( \pi(F)\right)_{\alpha \beta}= R \int_{X} F f_{\alpha,\beta} d\nu.
$$
Then the condition we are considering is that $\pi(\eta)\in {\cal H}(L^{k})$
should vanish, which is the same as saying that $\eta$ is orthogonal in
the $L^{2}$ sense to the finite dimensional subspace ${\rm Im}(\iota)\subset
C^{\infty}(X)$.
 
The practical question we now face is: how can we find the refined approximations
numerically, starting from the $\nu$ balanced metrics? As we have explained
above this is basically a linear problem. If we set $G=G_{0} + h$ then
we have a map
$$    V: h\mapsto \left(\eta_{\alpha \beta}(G_{0}+h)\right)  $$
from hermitian matrices $h$ to hermitian matrices $\eta_{\alpha \beta}$.
We want to find
a zero of $V$ and the standard procedure would be to invert the linearisation.
This linearisation is given by a $4$-index tensor $S_{\alpha \beta \gamma
\delta}$ with
$$ S_{\alpha \beta \gamma \delta}= \int_{X} f_{\alpha \beta} (f_{\gamma
\delta}+ \frac{1}{2} \Delta f_{\gamma
\delta}) d\nu. $$
 If we compute $S$ and invert the corresponding matrix we could define
an iterative procedure which ought to converge to the refined approximation.
An obstacle to carrying this through is that (even when reduced by symmetry)
the tensor $S$ has very many components so is relatively hard to compute
 in practice. Thus the author has not yet tried to implement this scheme
 but
 has used the following simpler procedure instead.
 Starting with an approximation $G_{0}$ we 
compute the error
matrix $E=\left( \eta_{\alpha \beta}\right)$ and simply set 
$$ G_{1}= G_{0}- \kappa E $$
where $\kappa$ is a suitable positive constant. In fact we compute  with the
inverse metrics $G_{0}^{-1}$ and, since $E$ will be small, use the approximation 
$$ G_{1}^{-1}= G_{0}^{-1}+ \kappa G_{0}^{-1} E G_{0}^{-1}. $$
Iterating this procedure yields a
sequence which appears to converges reasonably well, although slowly, to 
the refined approximation, see the examples and discussion
 in the next Section.
In Section 4 we give some further discussion and 
 theoretical justification for this procedure.

Of course
 the whole theory sketched here needs to be developed properly, and better
 methods found. But we hope it will yield systematic procedures for
 finding improved 
  numerical aproximations, starting with the balanced
 metrics.

\subsubsection{Non-zero cosmological constant}

       We now consider the problem of finding approximations to Kahler-Einstein
       metrics with non-zero scalar curvature. Of course this can be fitted
       into the constant scalar curvature theory described above, but there
       is also a natural variant of the Calabi-Yau construction above which
        probably yields a simpler approach (although the author has neither
        attempted to develop the theory of this nor studied substantial
        examples numerically). 
        
        We suppose that our positive line bundle $L$ is either the canonical
        bundle $K$ of $X$ or its dual $K^{-1}$. Write $p=\pm k$ in the
        two cases, so our space of sections is $H^{0}(X;K^{p})$. Given
        a metric  on $H^{0}(X;K^{p})$, let
             $s_{\alpha}$ be an orthornormal basis and set
      
        $$   \phi= \sum  s_{\alpha} \otimes \overline{s}_{\beta}. $$
        This is a section of the bundle $K^{p}\otimes \overline{K}^{p}$
        over $X$ which does not depend on the choice of orthonormal basis.
        Moreover if the sections of $K^{p}$ generate the fibers, which
        we can suppose is the case, $\phi$ does not vanish on $X$.
        Then $\phi^{1/p}$ is a well-defined volume form on $X$. Using this
        volume form in place of the Fubini-Study form we get the notion
        of a {\it canonically balanced} metric on $H^{0}(X;K^{p})$. Likewise,
        using this volume form we define another variant $T_{K}$ of the map $T$.
        Of course we hope that if a Kahler-Einstein metric exists then
        it is the limit of canonically balanced metrics, and that the 
        iterates of $T_{K}$ converge to the canonically balanced metrics.
        Furthermore, we can combine this discussion with that in (2.1.2), in the
        case when $X$ has continuous automorphisms, and we can hope to  find
        {\it Kahler-Ricci solitons} in appropriate cases. But we leave
        all of this as a programme for the future, except to give here
        another toy example. If we take invariant metrics on
         $\bC\bP^{1}$ and sections of ${\cal
        O}(6)= K^{-3}$ we get the sequence of iterates of $T_{K}$:
        
  $$        \begin{tabular}{|l||l|l|l|l|}\hline
        $r$& $a_{0}$ & $a_{1}$ &$a_{2}$& $a_{3}$\\ \hline
        0&.018&.5&4.5&54\\ \hline
        1& .0681&1.714&14.20&32.04\\ \hline
        2& .1860&3.120&15.50&26.39\\ \hline
        3&.3673&4.178&15.57&23.78\\ \hline
        4&.5598&4.900&15.43&22.23\\ \hline
        10&.9821&5.964&15.02&20.07\\ \hline
        18& .9998&6.000&15.00&20.00\\ \hline
        \end{tabular} $$
          
        The convergence is intermediate between the two previous cases,
        with the parameter $\sigma$  about .56.

\section{Study of a K3 surface}
\subsection{Geometry}
We now reach the heart of this article: the numerical study of a particular
K3 surface $S$. This is the double cover of the plane branched over the sextic curve
$x^{6}+ y^{6} + z^{6}=0$. Thus $S$ is defined by the equation $w^{2}= x^{6}+y^{6}+z^{6}$,
where $w$ is a point in the total space of the line bundle ${\cal O}(3)$
over $\bC\bP^{2}$. Most of the time we work in the affine piece of $S$
which is the subset of $\bC^{3}$, with co-ordinates $(x,y,w)$,
 cut out by the equation $w^{2}=x^{6}+y^{6}+1$. We fix the nowhere-vanishing
 holomorphic $2$-form $\theta$ on $S$ given in this affine piece by
 $$ \theta= \frac{dx dy}{w}= \frac{dx dy}{\sqrt{x^{6}+y^{6}+1}}. $$
This determines a volume form $\nu = \theta \overline{\theta}$.

The surface $S$ has many symmetries. These are generated by
\begin{itemize} \item The permutations of $x,y,z$
\item Multiplication of $x,y$ by sixth roots of unity.
\item The covering involution $w\mapsto -w$.
\item The antiholomorphic involution given by complex conjugation of all
coordinates.
\end{itemize}
Thus we get a symmetry group $\Gamma$ of order $6\times 6^{2} \times 2 \times 2=
864$, preserving the volume form $\nu$.

We now wish to evaluate the total mass
$$  {\rm Vol}(S,\nu)= \int_{S} \nu. $$
This is not strictly necessary for our main purpose, but gives a valuable
check on the accuracy of our numerical calculations in the next section.
We can  evaluate the integral by exploiting the fact that $S$ is an an elliptic surface,
of a very special kind. To see the elliptic fibration we consider the map
$ f:  [x,y,x]\mapsto [x^{3}, y^{3}, z^{3}]$ from $\bP^{2}$ to $\bP^{2}$.
This maps the sextic branch curve to the conic $X^{2}+Y^{2}+Z^{2}=0$. Thus
the covering $S\rightarrow \bP^{2}$ is the pull back by $f$ of the familiar
covering  of the quadric $Q$ over the plane, branched along the conic.
Now $Q= \bP^{1}\times \bP^{1}$ is fibred by lines (in two different ways).
These fibres are the preimages under the covering $Q\rightarrow \bP^{2}$
of the lines tangent to the conic. The preimages of these lines under
$f$ are cubic curves in the plane tangent to the sextic branch curve at
each intersection point, and these are the elliptic curves whose  lifts to
the double cover yield an elliptic fibration of $S$.

To see all of this more explicitly we work in affine coordinates. Let $C$
be the plane curve with equation $p^{3}+q^{3}+1=0$. Recall that the conic
$\lambda^{2}+\mu^{2}+1=0$ is parametrised by a rational variable $\tau$
with
\begin{equation} 
 \lambda=\frac{1}{2i}(\tau+\tau^{-1})\ \ , \ \ \mu=\frac{1}{2}(\tau-\tau^{-1}).
\label{eq:ratpar}\end{equation}
Now for a point $(p,q)$ in the curve $C$ and a complex parameter $\tau$
set
\begin{equation} x= \frac{p}{\lambda^{1/3}} \ \, \ \ y=\frac{q}{\mu^{1/3}}, \label{eq:cuber}\end{equation}
where $\lambda=\lambda(\tau), \mu=\mu(\tau)$ are given by Equation~\ref{eq:ratpar}. Of course
we have to deal with the cube roots in Equation~\ref{eq:cuber} , so initially we just regard
$\tau$ as varying in an appropriate open set $\Omega$ in $\bP^{1}$. 
A few lines of algebra show that if we set
$$  w= i \frac{ p^{3}-\lambda^{2}}{ \lambda \mu}, $$
then $x,y,w$ satisfy the equation $w^{2}=x^{6}+y^{6}+1$. In other words,
we have defined a holomorphic isomorphism  $F$ from $\Omega\times C\subset \bP^{1}
\times C$ to an open set in $S$. Further straightforward calculation shows
that
$$   F^{*}(\theta)= \frac{2^{2/3}}{3} \phi \psi $$
where $$ \phi=\frac{dp}{(1-p^{3})^{2/3}}\ \ , \ \ \psi= 
\frac{d\tau}{(\tau(\tau^{4}-1)^{1/3}}.
$$
The form $\phi$ is just the standard holomorphic $1$-form on the elliptic
curve $C$. The form $\psi$ is initially defined only an open set $\Omega$
in $\bP^{1}$. We could introduce a covering  $\Sigma\rightarrow
\bP^{1}$ such that $\psi$ lifts to a holomorphic form on $\Sigma$, but there
is no need to do this because the $2$-form $\psi \wedge \overline{\psi}$
is a well defined integrable form on $\bP^{1}$. Since we can cover a dense
open set in $S$ by taking $\Omega$ to be a cut plane, we see that
$$   {\rm Vol}(S,\nu) = \frac{2^{4/3}}{9} I  J, $$
where $$ I=i\int_{C} \phi\wedge \overline{\phi} \ \ , \ \ J=i\int_{\bP^{1}}
\psi \wedge \overline{\psi}. $$
Now let $W$ we the wedge-shaped region in $\bC$ defined by
$0<\arg(p)<\pi/3$. A dense open set in the curve $C$ is covered by 18 copies
of $W$ (6 rotations in $\bC$ each of which has 3 lifts to $C$). So
$$ I=18 \int_{W} \phi \wedge \overline{\phi}. $$
Now it is easy to see that the indefinite integral of $(1-p^{3})^{-2/3}$
maps $W$ to an equilateral triangle in $\bC$ with side length $L_{I}$ where
\begin{equation} L_{I}= \int_{0}^{1} \frac{dp}{(1-p^{3})^{2/3}}. \label{eq:integral}
\end{equation} 
It follows that
$$ I= 18 \frac{\sqrt{3}}{2} L_{I}^{2}. $$
Similarly let $Q$ be the domain in $\bC$ defined by the conditions
$\vert \tau \vert <1, 0<\arg(\tau)<\pi/2$. This is a fundamental domain
for an action of a group of order $8$ generated by  $\tau\mapsto \tau^{-1}$ and $\tau\mapsto
i\tau$ which preserves the form $\psi \wedge \overline{\psi}$. So
$$J=8 \int_{Q} \psi \wedge \overline{\psi}. $$
The indefinite integral of $\psi$ maps $Q$ to another equilateral triangle
with side length $L_{J}$ say, so
$$  J= 8 \frac{\sqrt{3}}{2} L_{J}^{2}, $$
where $$L_{J}=\int_{0}^{1} \frac{d\tau}{(\tau(1-\tau^{4}))^{1/3}}. $$
Finally, the elementary substitution,
\begin{equation}   p= \left(\frac{\tau^{2}-1}{\tau^{2}+1}\right)^{2/3}
\label{eq:subs}\end{equation}
shows that
$$  L_{J}= 3 \frac{4^{2/3}}{8} L_{I}. $$
(We can also interpret the substitution of Equation~\ref{eq:subs} as defining a covering map
from $\Sigma$ to $C$.) Putting everything together we obtain
$$  {\rm Vol}(S,\nu) = 27  L_{I}^{4}, $$
where $L_{I}$ is the one-dimensional integral in Equation~\ref{eq:integral}. Evaluating this numerically
one gets $L_{I}= 1.76664$, which yields
$$  {\rm Vol}(S,\nu)= 263.000, $$
accurate to about 5 or 6 significant figures. 
 
\subsubsection{Linear systems on $S$.}

We will study the general procedures of Section (2.2) for line bundles
${\cal O}(k)$ over $S$:  powers of the
 lift of the hyperplane bundle ${\cal O}(1)$ on $\bP^{2}$. We will consider
 three cases, when $k=3,6,9$. The symmetry group of $S$ acts on
 $H^{0}(S; {\cal O}(k))$ and we can restrict attention to invariant metrics.
 It is a straightforward matter to describe these. There is a natural holomorphic
 section $w$ of ${\cal O}(3)$ over $S$ and for $k\geq 3$ we have
 $$  H^{0}(S; {\cal O}(k)) = H^{0}(\bP^{2}, {\cal O}(k)) \oplus w H^{0}(\bP^{2};
 {\cal O}(k-3)). $$
 The two summands are eigenspaces of the action of the covering involution,
 so must be orthogonal for any invariant metric. The sections of ${\cal
 O}(j)$ over $\bP^{2}$ are represented, in our affine piece, by linear
 combinations of monomials 
 $x^{p} y^{q}$ for $p+q\leq j$. So we have a standard basis of the space
 labelled by the integer points in a triangle. We will draw this as a 
 right angled
 triangle but it is better to think of it as an equilateral triangle (with
 a hexagonal lattice) which makes the action of the  group of permutations
 of the projective coordinates 
 $x,y,z$ apparent. Now considering the action of multiplication by sixth
 roots of unity we see that we can only have a nontrivial inner product
 between monomials $x^{p} y^{q}$ and $x^{r} y^{s}$ if $p\equiv r$ and $q\equiv
 s$ modulo
 $6$. These inner products must be real numbers, due to the symmetry under
 the antiholomorphic involution of $S$. 
 
 Now for $k=3$ an invariant metric must be diagonal in our standard basis.
 The metric is specified by 4 real parameters $a_{I},a_{II},a_{III}, b_{I}$ as indicated.
 
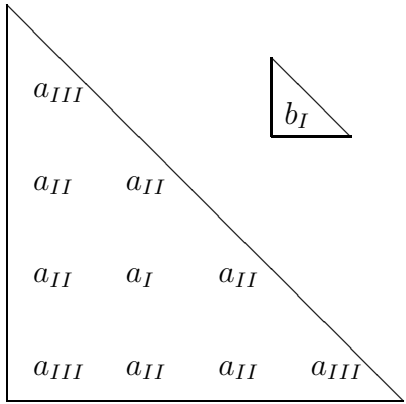
\begin{picture}(400,200)
 \put(0,0){\line(1,0){150}}
 \put(0,0){\line(0,1){150}}
 \put(150,0){\line(-1,1){150}}
 \put(10,10){$a_{III}$}
 \put(45,10){$a_{II}$}
 \put(80,10){$a_{II}$}
 \put(115,10){$a_{III}$}
 \put(10,45){$a_{II}$}
 \put(10,80){$a_{II}$}
 \put(10,115){$a_{III}$}
 \put(45,45){$a_{I}$}
 \put(80,45){$a_{II}$}
 \put(45,80){$a_{II}$}
 \put(100,100){\line(1,0){30}}
 \put(100,100){\line(0,1){30}}
 \put(130,100){\line(-1,1){30}}
 \put(105,105){$b_{I}$}
 \end{picture}

 In other words, the function $D$ associated to the metric is given, in
 our affine coordinates, by
 $$D= a_{I} \vert xy\vert^2 +a_{II} (\vert x\vert^2+\vert x\vert^4+\vert y \vert^2+\vert
 y\vert^4+ \vert x^2 y\vert^2 + \vert y^2 x\vert^2) + a_{III} (\vert x\vert^6
 +\vert y\vert ^6+ 1 ) + b_{I} \vert
 w\vert^2 $$
 
 For $k=6$ we again have diagonal elements filling up two triangles, invariant
 under the permutations of $x,y,z$. This gives us 10 parameters $a_{I},a_{II},\dots
a_{VII}$,\ $ b_{I},b_{II},b_{III}$ according
 to the scheme indicated:

 
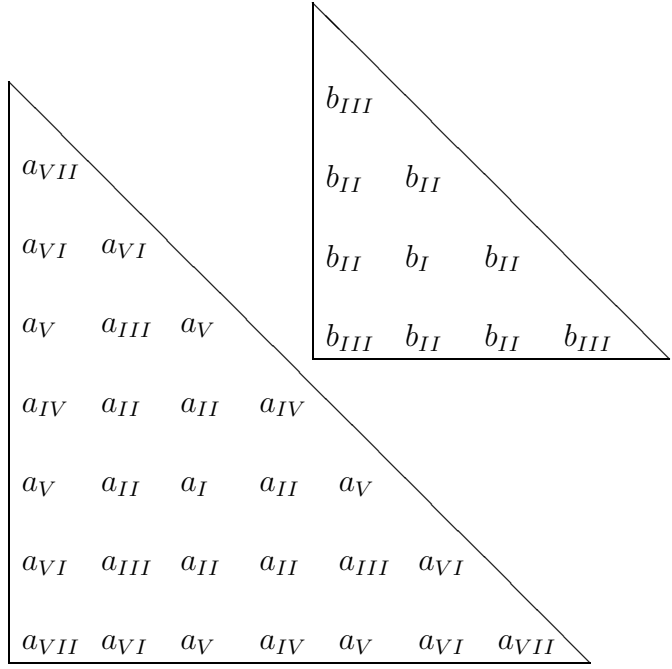
\begin{picture}(500,300)
 \put(0,0){\line(1,0){220}}
 \put(0,0){\line(0,1){220}}
 \put(220,0){\line(-1,1){220}}
 \put(65,65){$a_{I}$}
 \put(5,5){$a_{VII}$}
 \put(35,5){$a_{VI}$}
 \put(65,5){$a_{V}$}
 \put(95,5){$a_{IV}$}
 \put(125,5){$a_{V}$}
 \put(155,5){$a_{VI}$}
 \put(185,5){$a_{VII}$}
 \put(5,35){$a_{VI}$}
 \put(5,65){$a_{V}$}
 \put(5,95){$a_{IV}$}
 \put(5,125){$a_{V}$}
 \put(5,155){$a_{VI}$}
 \put(5,185){$a_{VII}$}
 \put(35,155){$a_{VI}$}
 \put(65,125){$a_{V}$}
 \put(95,95){$a_{IV}$}
 \put(125,65){$a_{V}$}
 \put(155,35){$a_{VI}$}
 \put(35,35){$a_{III}$}
 \put(125,35){$a_{III}$}
 \put(35,125){$a_{III}$}
 \put(35,65){$a_{II}$}
 \put(65,35){$a_{II}$}
 \put(35,95){$a_{II}$}
 \put(95,35){$a_{II}$}
 \put(65,95){$a_{II}$}
 \put(95,65){$a_{II}$}
\put(115,115){\line(1,0){135}}
 \put(115,115){\line(0,1){135}}
 \put(250,115){\line(-1,1){135}}
 \put(120,120){$b_{III}$}
 \put(150,120){$b_{II}$}
 \put(180,120){$b_{II}$}
 \put(210,120){$b_{III}$}
 \put(120,150){$b_{II}$}
 \put(120,180){$b_{II}$}
 \put(120,210){$b_{III}$}
 \put(150,150){$b_{I}$}
 \put(180,150){$b_{II}$}
 \put(150,180){$b_{II}$}
 \end{picture}
 
 We also have three non-trivial off-diagonal terms corresponding to the inner products
 between $1, x^6, y^6$. These must all be equal, by symmetry, so we get
 one further real parameter $C$. In other words the matrix of our inner
 product contains a $3\times 3$ block
 $$ \left(\begin{array}{ccc}   a_{VII}&C&C\\C&a_{VII}&C\\C&C&a_{VII}
 \end{array}\right).$$
  
 (We could diagonalise this by using a different basis, but to fit in with
 the notation below we will not do so.)
                   
In other words, the function $D$ is given by
$$  D= \sum a_{p,q} \vert x\vert^{2p} \vert y \vert^{2q} + C( x^6 \overline
{y}^{6} +\overline{x}^{6}y^{6} + x^{6}+\overline{x}^{6}+y^{6}+\overline{y}^{6}) + \vert w\vert^{2} \sum b_{p,q} \vert x\vert^{2p}
\vert y\vert^{2q}, $$
where the coefficients $a_{p,q}, b_{p,q}$ are given by the parameters
$a_{I},\dots, a_{VII} $, $b_{I}, \dots, b_{III}$ in the manner  indicated above.

Finally we consider the case $k=9$. The diagonal elements in the big triangle
are specified by twelve parameters $a_{I},\dots , a_{XII}$ in the manner indicated
(here we only draw a piece of the big triangle, the remainder follows by
symmetry).

\begin{picture}(300,150)
\put(0,0){\line(1,0){250}}
\put(0,0){\line(0,1){150}}
\put(5,5){$a_{XII}$}
\put(35,5){$a_{XI}$}
\put(65,5){$a_{X}$}
\put(95,5){$a_{IX}$}
\put(125,5){$a_{VIII}$}
\put(155,5){$a_{VIII}$}
\put(95,95){$a_{I}$}
\put(95,65){$a_{II}$}
\put(65,65){$a_{III}$}
\put(125,35){$a_{IV}$}
\put(95,35){$a_{V}$}
\put(65,35){$a_{VI}$}
\put(35,35){$a_{VII}$}
\end{picture}

There are now more allowable off-diagonal elements, specified by  6 
independent parameters $C_{1}, \dots, C_{6}$. These fit into $3\times 3$ blocks
in the matrix of the metric of the following forms:

$$ \left(\begin{array}{ccc} a_{VI}&C_{2}&C_{1}\\ C_{2}&a_{X}& C_{3}\\ C_{1}&C_{3}&a_{XI}\end{array}\right)\
\ , \left(\begin{array}{ccc} a_{XII} &C_{4}&C_{4}\\C_{4}&a_{IX}&C_{5}\\
C_{4}&C_{5}&a_{IX}\end{array}\right)\ , \ 
\left(\begin{array}{ccc} a_{VII} &C_{6}&C_{6}\\C_{6}&a_{VII}&C_{6}\\C_{6}&C_{6}&a_{VII}\end{array}\right)
. $$
The discussion of the small triangle repeats that in the $k=6$ case above:
we have 8 parameters $b_{I},\dots ,b_{VII}, C'$.

To sum up,
we have now specified how invariant metrics on $H^{0}(S; {\cal O}(k))$
are determined by 4,11 and 26 real parameters respectively  in the three cases $k=3,6,9$.
Notice that these are vast reductions on the dimensions 
of the full spaces of metrics ${\rm dim}\ H^{0}(S;{\cal O}(k))^{2}$, which
are 121, 1444,6889 respectively.

\subsection{Numerical results}
\subsubsection{Numerical volume}
The numerical implementation of the algorithm of (2.2) for finding  a balanced
metric is completely specified by the choice of an approximation
$\nu^{*}=\sum \nu_{i} \delta_{z_{i}}$ for the given volume form $\nu$ on
$S$. That is, what we compute is the sequence $T_{\nu^{*}}^{r}(H_{0})$.
Our numerical results are obtained using a family of such approximations which depend
on four integer parameters $n_{x}, n_{p}, n_{u},n_{w}$: the larger these
parameters are the  more points $z_{i}$ are used, which should give a better approximation,
at the expense of extra computation. More details are given in the next
section, suffice it to say  here that $n_{x}, n_{p}$ determine the contribution
from one chart and $n_{u}, n_{w}$ from another. We should also emphasise
that our approximation exploits the invariance of the functions we need
to integrate under the group of symmetries of $S$.

As a first test of our approximations we compute the volume ${\rm Vol}(S,\nu_{0})$ of $S$
in the measure $\nu$. We obtain then a family of
 approximations $V(n_{x}, n_{p}, n_{u}, n_{w})$ with
 $$ V(n_{x}, n_{p}, n_{u}, n_{w})= V_{1}(n_{x}, n_{p}) + V_{2}(n_{u}, n_{w}).
 $$
 We find:

$$ \begin{tabular}{|l||l|}
 \hline $V_{1}(10,10)=265.84$& $V_{2}(10,10)=6.5557$\\ \hline 
 $V_{1}(20,20)=256.61$&$V_{2}(14,10)=6.3227$\\ \hline 
 $V_{1}(30,30)=256.88$&$V_{2}(20,20)=6.2575$\\ \hline 
 $V_{1}(40,40)=256.70$&$V_{2}(24,20)=6.2536$\\ \hline
 \end{tabular} $$
 
So, for example, $V(40,40,24,20)=256.70+6.25=262.95$ and
\linebreak $V(20,20,14,10)=262.93$. Recall that the theoretical analysis gave
${\rm Vol}(S,\nu)=263.000$. These figures suggest that, for values of the parameters
similar to those above, we can integrate reasonably smooth functions on
$S$ with an accuracy of about 3 or 4 significant figures. 
   
\subsubsection{The case $k=3$}

We now implement our \lq\lq $T_{\nu}$-algorithm'' to find the $\nu$-balanced
metric for the line bundle ${\cal O}(3)$ on $S$. We choose, arbitrarily,
the initial metric specified by the four parameters $(1,1,1,1)$. Then we
obtain, using the approximating measure $\nu^{*}(20,20,14,10)$:

$$\begin{tabular}{|l||l|l|l|l|}
\hline
r& $a_{I}$&$a_{II}$&$a_{III}$&$b_{I}$ \\ \hline
0&1&1&1&1\\ \hline
1&12.68&8.758&5.257&2.722\\ \hline
2&13.20&8.807&4.987&2.439\\ \hline
3& 13.25&8.811&4.959&2.414\\ \hline
4&13.26&8.812&4.956&2.412\\ \hline
5&13.26&8.812&4.956&2.412\\ \hline
\end{tabular}$$
(Recall that everything is preserved by rescaling the metrics.
Here, and in the similar results to follow, we normalise the metrics
for $r\geq 1$, up to scale, in an arbitrary way.) Rounding off to 4 SF, the iteration reaches a fixed point after four steps.
Taking a finer approximation $\nu^{*}(30,30,20,14)$ makes hardly any change:
to 4SF, the fixed point is now 13.26,8.816,4.956,2.415.
So we have some  confidence that these last values of the parameters give
the balanced metric on $H^{0}(S;{\cal O}(3))$ to high  accuracy.

We now examine the Fubini study volume form $\mu$ on $S$ determined by
this balanced metric. (To fix constants, we take the Kahler form in the
class $2\pi c_{1}({\cal O}(3))$.) That is, for given parameters $n_{x}, n_{p}, n_{u},
n_{w}$ we compute the ratio $ \mu/\nu$ at each of the points in the support
of $\nu^{*}(n_{x}, n_{p}, n_{u}, n_{w})$ (see Section (3.4.2) for more details
of the calculation). Now we have another test of our numerical methods
by integrating $\mu/\nu$ with respect to the measure $\nu^{*}(n_{x}, n_{p},
n_{u}, n_{w})$. By Chern-Weil theory 
$$ \int_{S} d\mu = 9.2.4.\pi^{2} = 710.61. $$
so the ratio of the total volumes with respect to $\mu$ and $\nu$
is $710.61/263=2.7019$.
Taking parameters $[28,28,20,14]$ we get
$$ \int_{S} \frac{\mu}{\nu} d\nu^{*}= 2.7023 \int_{S} d\nu^{*}  $$
which is again a fair agreement. Set $\eta=  (2.7019)^{-1} \mu/\nu$, so the
mean  value of $\eta$ with respect to $\nu$ is $1$ and the deviation
of $\eta$ from the constant function $1$ is a measure of the difference
between the balanced metric and the Calabi-Yau metric on $S$, normalised
by the appropriate scale factor. We compute
\begin{itemize}
\item The maximum value of $\eta$,
\item The minimum value of $\eta$,
\item The mean value of $\vert \eta-1\vert$ with respect to $\nu$, 
\item The distribution function of $\eta$ with respect to $\nu$.
\end{itemize}
More precisely, of course,  we compute the maximum and minimum values over the support
of a $\nu^{*}(n_{x}, n_{p}, n_{u}, n_{w})$, for appropriate parameter values,
and we compute the mean and distribution function with respect to $\nu^{*}$.
The result (with parameters [28,28,20,14]) is 
$$  \max \eta=1.496,\  \ \min \eta= .2501,\ \  
{\rm Mean}(\vert \eta-1\vert)= .262$$

For the (approximate) distribution function we give in the first row of
the table a collection of ranges and in the second row the percentages
of the $\nu$-volume of $S$ where $\eta$ takes values in the given range.

$$\begin{tabular}{|l|l|l|l|l|l|l|l|l|}\hline
.25--.4&.4--.55&.55--.7&.7-.85&.85--1&1--1.15&1.15--1.3&1.3-1.45&1.45--\\ 
\hline
 3.4  &6.9&9.6    &11.7    &15.3  &17.3    &15.0      &
15.8 &5.1   \\ \hline
\end{tabular}$$

We see that $\mu$ is not really a good approximation to $2.7019 \nu$, 
as we would expect
with this very low value of $k$. So we do not proceed to try to
find a  \lq\lq refined approximation'' in this case.

\subsubsection{The case $k=6$}
We carry through the same procedure as in the previous subsection. With
parameters $[20,20,14,10]$ the following table shows the convergence to
the balanced metric.

$$\begin{tabular}{|l||l|l|l|l|l|l|l|l|l|l|l|}
\hline
r&$a_{I}$&$a_{II}$&$a_{III}$&$a_{IV}$&$a_{V}$&$a_{VI}$&$a_{VII}$&$b_{I}$&$b_{II}$&$b_{III}$&$C$\\ \hline
0&1&1 &1  &1 &1&1 &1  &1 &1  &1   & 0\\ \hline
1&51.7&43.3 &33.1  &23.9 &21.8&15.9 & 7.55 &21.8 &14.1 &7.37   &.915\\
 \hline
 2&56.4&45.9&33.6&23.2&20.8&14.3&5.85&20.3&12.2&5.66&.509\\
 \hline
 3&57.6&46.5&33.7&23.0&20.6&13.9&5.47&19.9&11.7&5.28&.432\\
 \hline 4&57.6&46.6&33.7&23.0&20.6&13.9&5.46&19.9&11.7&5.27&.430\\
 \hline $\infty$&57.6&46.6&33.7&23.0&20.6&13.9&5.45&19.9&11.7&5.26&.429\\
 \hline
\end{tabular}
$$
On the author's PC each step takea about 5 minutes.
Analysis of the convergence suggests that the parameter $\sigma$ is about
.22.
Increasing the numerical integration parameters to $[30,30,20,16]$
gives a fixed point, to 4 significant figures, 

$$\begin{tabular}{|l|l|l|l|l|l|l|l|l|l|l|}\hline
57.70&46.61&33.67&23.01&20.56&13.86&5.453& 19.95&11.70&5.267&  .4324 \\
\hline 
\end{tabular}
$$
We take these parameters as our numerical balanced metric.

The volume form of the balanced metric yields a function $\eta$ with:
$$ \max \eta= 1.065\ ,\  \min \eta= .679\ ,\    {\rm Mean}(\vert \eta-1\vert)= .058.
$$
The distribution function of $\eta$ is given by the following table, in
which again the first row gives the range and the second gives the percentage
of the total volume lying within the range.

\begin{tabular}{|l|l|l|l|l|l|l|l|l|}
 \hline
-.7&.7-.75&.75-.8&.8-.85&.85-.9&.9-.95&.95-1&1-1.05&1.05-\\
\hline
.07&1.2   &2.2   &3.3   &4.9   &7.9   &14.0 &35.0  &31.4 \\ \hline
\end{tabular}

This is a much better approximation to the Calabi-Yau metric. The function
$\eta$ is close to $1$ over most of the manifold; further exploration shows
 the set where it deviates
substantially from $1$ is a neighbourhood of the branch curve of the double
cover, on which $\eta$ is small. 

We now move on to search for a refined approximation, using the algorithm
described in (2.2.1) above. Thus we compute the tensor $\pi(\eta)\in {\cal
H}(L^{k})$
as a  measure of the size of the error term. If we work with an orthonormal
basis of sections then the individual 
matrix entries $\eta_{\alpha \beta}$ can be interpreted as analogues of
Fourier coefficients of $\eta$. For simplicity we work instead with an
approximately orthonormal basis given by rescaling the standard monomials,
which we justify by the fact that the off-diagonal term $C$ is relatively
small and we are working close to the balanced metric. Thus for any metric we define
$$ \eta(pq)= R a_{pq} \int_{S} (\eta-1) \frac{\vert x\vert^{2p} \vert y \vert^{2q}}{D}
d\nu, $$
$$ \overline{\eta}(pq)= R
b_{p,q} \int_{S} (\eta-1) 
\frac{ \vert x\vert^{2p} \vert y^{2p} \vert w \vert^{2}}{D}
$$ and $$ \eta_{C}= R a_{0,0} \int_{S} (\eta-1) \frac{{\rm Re}(x^{6})}{D}
d\nu_{0} $$
where $R= {\rm dim H^{0}(L^{k})}/{\rm Vol}(S,\nu_{0})= 38/263$. The terms
$\eta(pq), \overline{\eta}(pq)$ have the same symmetries as the coefficients
so we can extract 11 different terms $\eta_{I}, \dots, \eta_{VII}, \overline{\eta}_{I},\dots,
\overline{\eta}_{III}, \eta_{C}$ corresponding to the coefficients
$a_{I}\dots,a_{VII}$,\ $ b_{I},\dots, b_{III}, C$. As described in Section
2, the diagonal terms $\eta_{I},\dots, \overline{\eta}_{III}$ can be interpreted
approximately as weighted averages of the error $\eta-1$. 
(The approximations we are making here only involve the interpretetation
of the data, not the actual algorithms.)

 The next two tables display the first 5 steps of the refining procedure
with the  parameter $\kappa$ equal to 2.5. The sequence of metrics is:

$$\begin{tabular}{|l||l|l|l|l|l|l|l|l|l|l|l|}
\hline
r&$a_{I}$&$a_{II}$&$a_{III}$&$a_{IV}$&$a_{V}$&$a_{VI}$&$a_{VII}$&$b_{I}$&$b_{II}$&$b_{III}$&$C$\\ \hline
0&57.70&46.61 &33.67  &23.01 &20.56&13.86 &5.452  &19.95 &11.70  &5.267&  .4326\\ \hline
1&53.63&44.15 &33.47  &23.92 &21.62&15.08 & 6.380 &21.73 &13.19 &6.124   &.6059\\
 \hline
 2&51.62&43.72&33.51&24.18&21.88&15.26&6.394&22.52&13.55&6.136&.6100\\
 \hline
 3&51.01&43.50&33.59&24.34&22.03&15.34&6.347&23.17&13.74&6.096&.5949\\
 \hline 4&50.57&43.35&33.62&24.44&22.13&15.39&6.310&23.61&13.87&6.066&.5785\\
 \hline 5 &50.22&43.24&33.74&24.51&22.19&15.43&6.284&24.03&13.95&6.047&.5621\\
 \hline
\end{tabular}$$
The corresponding sequence of $\eta$-coefficients, multiplied by $10^{3}$:
 
$$\begin{tabular}{|l||l|l|l|l|l|l|l|l|l|l|l|}
\hline
r&$\eta_{I}$&$\eta_{II}$&$\eta_{III}$&$\eta_{IV}$&$\eta_{V}$&$\eta_{VI}
$&$\eta_{VII}$&$\overline{\eta}_{I}$&$\overline{\eta}_{II}$&$\overline{\eta}_{III}$&$\eta_{C}$\\ 
\hline
0&-56.1&-42.9 &-25.2  &-8.0 &-3.6&10.1 &32.0 &10.7 &25.1  &38.3&  3.0\\ \hline
1&-10.7&-6.9 &-2.6  &1.2 &1.6&1.7 & -4.8 &11.2 &7.6 &-2.2   &.3\\
 \hline
 2&-5.9&-3.3&-.4&1.3&1.4&.8&-5.1&9.3&4.4&-3.9&-.2\\
 \hline
 3&-4.4&-2.3&.03&.8&.9&.5&-3.3&7.7&2.8&-2.8&-.4\\
 \hline 4&-3.4&-1.6&.2&.5&.6&.3&-1.8&6.4&1.7&-1.9&-.6\\
 \hline 5 &-2.8&-1.3&.4&.3&.3&.2&-.8&5.4&.9&-1.3&-.7\\
 \hline
 \end{tabular}$$
 
The size of the error terms behaves as follows:

\begin{tabular}{|l||l|l|l|}\hline
r& Max. &Min.& Mean  Error(\%) \\ \hline
0& 1.065    &.679  &5.80\%    \\ \hline
1& 1.046    &.804  &2.61\%    \\ \hline
2& 1.043    &.830  &2.21\%    \\ \hline
3& 1.041    &.843  &1.94\%   \\ \hline
4&1.039&.853&1.79\% \\ \hline
5& 1.05&.860&1.75\% \\ \hline
\end{tabular}

We see that for the first four steps the error term does 
 decrease, according to any of the
measures above. The initial rate of decrease is large but this soon slows
down. At the fifth step the maximum value increases, although the mean
and minimum improve.  The most obvious phenomenon is
that the coefficient $a_{I}$ is decreased, along with  the \lq\lq nearby'' co-efficients in the
big triangle, and the co-efficient $b_{I}$ is increased, along with the nearby
coefficients in the small triangle. This has the effect of increasing the
volume form $\mu$ near the branch curve, and so compensating for the deviation
in the balanced metric. We repeat the process many times, with various
values of $\kappa$. After a while
the maximum value decreases again, and for example we achieve at an intermediate
stage the
 following  metric, which we call $\omega_{6}'$: 
$$\begin{tabular}{|l|l|l|l|l|l|l|l|l|l|l|}\hline
46.61&42.69&34.94&25.04&22.56&15.34&6.301&28.95&13.79&6.130&.2854 \\ \hline
\end{tabular}$$
This has 
$$ {\rm Max}\ \eta = 1.031 \ ,\ {\rm Min.}\ \eta= .898\ ,\  {\rm Mean}(\vert \eta -1\vert)=
.017  , $$
and $\eta$-coefficients, multiplied by $10^{3}$:
$$\begin{tabular}
{|l|l|l|l|l|l|l|l|l|l|l|}\hline
-.57&-.11&.1&.11&.04&-.3&.25&.9&-.2&.1&-.2
\\ \hline
\end{tabular}$$
The distribution function is:
$$\begin{tabular}
{|l|l|l|l|l|l|l|l|l|} \hline
--.925&.925-.94&.94-.955&.955-.97&.97-.985&.985-1&1-1.015&1.015-1.03&1.03-\\
\hline
.32&.45&1.0&1.4&7.9&38.8&33.4&16.7&.004 \\ \hline
\end{tabular} $$

We contine the process still further. The rate of change becomes extremely
slow, and we finally give up at a metric we will call $\omega_{6}''$ defined
by metric coefficients
$$ \begin{tabular}{|l|l|l|l|l|l|l|l|l|l|l|}\hline
43.34&42.77&35.84&25.31&22.86&14.88&6.362&30.19&13.78&6.207&.1341\\
\hline
\end{tabular}$$
Which has
$$ {\rm Max}\ \eta=1.050\ ,\  {\rm Min}\ \eta=.909\ ,\  {\rm Mean}(\vert \eta-1\vert)=.0111$$
and $\eta$-coefficients, multiplied by $10^{3}$:
$$\begin{tabular}
{|l|l|l|l|l|l|l|l|l|l|l|}\hline
-.21&.011&.063&-.026&.020&-.070&-.004&-.063&.024&.005&-.047\\
\hline \end{tabular}$$
We take this metric $\omega''_{6}$ as our best estimate for a refined 
approximation. Now there are several points to make. First, the
 $\eta$-coefficients for $\omega_{6}''$ are much smaller than those
 for the
balanced metric $\omega_{6}$, by a factor of more than $100$.
 However
the very slow movement under our primitive algorithm means that it is not completely
clear that there is a genuine, exact,  refined approximation close to
$\omega_{6}''$. Second, the metric $\omega''_{6}$ is probably not significantly
 better
as an approximation to the Calabi-Yau metric than the intermediate metric
$\omega_{6}'$, and could well be worse. So our overall conclusion is that
the ideas we developed in (2.2.1) are effective, in this case, in generating
some much better approximations than the balanced metric (such as $\omega'_{6}$
or even the fourth step in the process), although a more exact and sophisticated
 analysis of the whole situation is required.  

\subsection{The case $k=9$.}

The numerical results in this case follow very much the same pattern as
for $k=6$ above, so we will be brief. The crucial difference, of course, is that the error terms
are smaller. Recall that we now have $26$ parameters which we display in
four-row form as

$$\begin{tabular}{|l|l| l|l|l|l|l|} \hline
 $a_{I}$&$a_{II}$&$a_{III}$&$a_{IV}$&$a_{V}$&$a_{VI}$&$a_{VII}$\\ \hline
 $a_{VIII}$&$a_{IX}$&$a_{X}$&$a_{XI}$&$a_{XII}$&\ &\ \\ \hline
 $b_{I}$&$b_{II}$&$b_{III}$&$b_{IV}$&$b_{V}$&$b_{VI}$&$b_{VII}$\\ \hline
  $C_{1}$&$C_{2}$&$C_{3}$&$C_{4}$&$C_{5}$&$C_{6}$&$C'$\\ \hline
\end{tabular}$$

We find the balanced metric is:
$$
\begin{tabular}{|l|l|l|l|l|l|l|} \hline
 992.0& 858.9& 682.7& 572.8& 525.7& 401.8 &227.6\\ \hline
 200.8& 176.9& 130.3& 69.46& 19.60&\ &\ \\ \hline
 460.1& 346.2& 218.2& 147.8& 123.8& 68.75& 19.57 \\ \hline
18.45&19.26&17.82&7.275&8.245&33.21&.1903 \\ \hline
\end{tabular}
$$
The convergence parameter $\sigma$ is about .33. The balanced metric has
 
$$ \max \eta= 1.024 \ , 
\ \min \eta=.833\ ,\  {\rm Mean}(\vert \eta-1\vert)= .025. $$
The largest modulus of the $\eta$-coefficients is about $34\times 10^{-3}$.
 
The distribution function is

$$\begin{tabular}{|l|l|l|l|l|l|l|l|l|}
 \hline
-.88&.88-.9&.9-.92&.92-.94&.94-.96&.96-.98&.98-1&1-1.02&1.02-\\
\hline
1.3 &1.7&2.1&3.1&4.2&6.6&12.1&29.9&39.1 \\ \hline
\end{tabular}$$

Thus over about 70\% of the volume, $\eta$ lies in the range $1-1.025$.  
We carry out the refining procedure as before. The process seems to work
better than for $k=6$, in that the error decreases steadily and our best
estimate for the refined approximation  yields the best approximation
to the Calabi-Yau metric. This approximation, $\omega'_{9}$ is
$$ \begin{tabular}{|l|l|l|l|l|l|l|} \hline
 875.6& 798.5& 665.8& 584.5& 539.3& 419.5 &239.8\\ \hline
 214.3& 190.4& 140.3& 76.43& 21.76&\ &\ \\ \hline
 540.5& 386.8& 233.3& 161.8& 134.9& 75.53& 21.76 \\ \hline
19.89&25.89&19.57&8.406&14.44&31.69&.5101 \\ \hline
\end{tabular}
$$

The $\eta$-coefficients (times $10^{3}$) are  
$$\begin{tabular}{|l|l|l|l|l|l|l|}
\hline
-.09&-.04&-.001&-.02&-.02&-.03&-.05\\
\hline
-.03&-.04&-.03&.005&-.06&\ \\
\hline
.16&-.002&-.08&-.02&-.02&.01&-.04\\
\hline
-.01&.02&.03&.03&.04&-.02&.03\\ \hline
\end{tabular}$$

The error has
$$ \max \eta= 1.009\  ,\  \min \eta= .972\ ,\  {\rm Mean}(\vert \eta-1\vert)= .0022, $$
and distribution function:

$$\begin{tabular}{|l|l|l|l|l|l|l|l|}
\hline
-.975&.975-.98&.98-.985&.985-.99&.99-.995&.995-1&1-1.005&1.005-\\ \hline
  .01&.1&.3&.9&1.5&46.6&46.1&4.5\\ \hline
\end{tabular}$$

So over more than 90\% of the manifold $\eta$ is within .5\% of 1.
Notice that in passing from $k=6$ to $k=9$, we reduce the mean error in
the balanced metrics
by a factor of about $2$, which is roughly $(9/6)^{2}$. On the other hand
we reduce the mean error in the refined approximations by a factor of about
$5$. This gives some small support for the idea that the refined approximations
should yield rapid convergence, in $k$, to Calabi-Yau metrics.

This concludes our account of the main numerical results of this article.
Notice that the function $\vert \theta\vert_{\omega}$ discussed in the
Introduction is  essentially $\eta^{-1/2}$, so the deviation from $1$ will
be about half that of $\eta$.

\subsection{Computational details}
\subsubsection{Charts}
To integrate numerically over $S$ we need an appropriate system of coordinate
charts and choice of good charts is an interesting geometrical question. Suppose $(\xi,
\eta)$ are local holomorphic coordinates for some chart in an atlas, i.e with
$(\xi, \eta)$ lying in some bounded  domain $D\subset \bC^{2}$ and with the holomorphic
$2$-form given by $\theta = F(\xi,\eta) d\xi d\eta$. Then we want to realise
three desireable properties:
\begin{itemize}
\item The oscillation of the function $F$ is not too large on $D$.
\item It is easy to recognize if a point $(\xi,\eta)$ lies in $D$.
\item The overlaps of this chart with the other charts in the atlas are
not too small.
\end{itemize}
The last property is needed because our integration procedure will involve
a partition of unity subordinate to the cover, so we need a smooth cut-off
function $\psi_{D}$ supported in $D$ and equal to $1$ outside the region
corresponding to the other charts. The contribution from $D$ to the integral of
a function $f$ on $S$ will have the shape
$$ \int_{D} f \psi_{D} \vert F\vert^{2} d\xi d\overline{\xi} d\eta d{\overline
\eta}. $$
If the overlaps are small the derivative
of $\psi_{D}$ will be large and the numerical integration will not be accurate.

As a first step we restrict to the open set $U$ in $S$ defined by the conditions
$$  \vert x\vert <1.2 \vert z \vert\ ,\ \vert y\vert < 1.2 \vert z \vert,
$$
in terms of the homogeneous co-ordinates $(x,y,z)$ on $\bC\bP^{2}$. Clearly
$S$ is covered by $U$ and its two images under the permutation action,
so to integrate a $\Gamma$-invariant  function it suffices  
to work in $U$. Then we can pass to affine coordinates, where $z=1$ and
$U$ is represented as the polydisc $\vert x\vert, \vert y\vert<1.2$. 
Now let $U_{1}$ be the complement in $U$ of the product of the disc
$\vert y \vert<1.2$ with suitable small neighbourhoods of the points where
$x^{6}=-1$, chosen so that the roots of $(1+x^{6})$ can be defined on
$U_{1}$. Then
$$   w^{2}= (1+x^{6}) \left( 1 +\left(\frac{y}{(1+x^{6})^{1/6}}\right)^{6}\right).
$$
Set $q=\frac{w}{\sqrt{1+x^{6}}}$ so that
$$ q^{2}-1 = \frac{y^{6}}{1+x^{6}}. $$
Now let $U_{2}\subset U_{1}$ be the region defined by the condition ${\rm
Re}(q)>-1/2$.Then $U_{1}$ is covered by $U_{2}$ and its image under the covering involution,
(which takes $q$ to $-q$).  On $U_{2}$ we write $q=1+p^{6}$ so
 \begin{equation} y= p(2+p^{6})^{1/6} (1+x^{6})^{1/6}, \label{eq:ydef}\end{equation}
where the root $(2+p^{6})^{1/6}$ is well-defined on $U_{2}$ since we have
cut out the zero where $q=-1$ and $2+p^{6}=0$. The upshot is that we have
an open set $U_{2}\subset S$ on which we can take $x$ and $p$ as complex
co-odinates, with $y$ given by Equation~\ref{eq:ydef} and 
$$  w= (1+p^{6}) \sqrt{1+x^{6}}. $$
This gives one chart, which we call the \lq\lq big chart'' in $S$. 
 The holomorphic form in these coordinates is 
$$ \frac{2}{ (2+p^{6})^{5/6} (1+ x^{6})^{1/3}}dp dx.  $$
Interchanging $x$ and $y$ takes $U_{1}$ to another open set $U'_{1}\subset
U$ say. The union $U_{1} \cup U'_{1}$ covers all of $U$ save for small
neighbourhoods of the points where $x^{6}=-1$ and $y^{6}=-1$. We want to
define coordinates on an open set $V$ in $S$
 covering a neighbourhood of the points $x=y=e^{i\pi/3}$. Then $U$ will
 be covered by $U_{1}, U'_{1}$ and the 36 images of $V$ under the action
 of multiplication by sixth roots of unity. To integrate a $\Gamma$-invariant
 function it will suffice to work in $U_{2}$ and $V$, provided we use suitable
 invariant cut-off functions and take due account of multiple counting
 by the symmetry group. 
 
 On $V$ we set $u=x^{6}-y^{6}$ and take $u$ and $w$ as local coordinates so
 $$   x= \left( \frac{w^{2}+u-1}{2}\right)^{1/6}\ , \ y=\left( \frac{w^{2}-u-1}{2})\right)^{1/6}.
 $$
 This gives our other chart, the \lq\lq small chart'', on $S$. 
  in which the holomorphic form is 
 $$  \frac{du dw}{36 x^{5} y^{5}}. $$
The tension set up by the requirement discussed above for the two charts can now
be seen as follows. In the small chart we need $\vert x\vert$ not to 
be too small and  in the large chart we need $1+x^{6}$
not to be too small, but  the two charts must have substantial overlap.
We balance these requirements by defining the domain of the big chart to
be the region where ${\rm Re}(x^{6})>-(.9)^{6}=-.531 $ and the domain of the small
chart to be where ${\rm Re}(x^{6})<-(.7)^{6}=-.117$.
 Then in the big chart $\vert 1+x^{6}\vert>.47$, and  in
the small chart $\vert x\vert>.7$. The  cut-off function has
 derivative (with respect
to $x$) about $2(.9-.7)^{-1}=10$.
       
       To integrate numerically in the big chart we use hexagonal lattices
       in the complex $x$ and $p$ variables. These allow us to take exact
       account of the residual $\bZ/6\times \bZ/6$ symmetry (that is the
       integrands are functions of $x^{6}$ and $p^{6}$) saving a factor
       of nearly 36 in the calculation. We nearly double this factor using
       the further symmetry under complex conjugation.   The lattice spacings
       are proportional to $n_{x}^{-1}, n_{p}^{-1}$ where $n_{x}, n_{p}$
       are integer parameters mentioned in (3.2.1). To integrate numeriaclly
       in the small chart we use square lattices in the $u$ and $w$
        variables, taking advantage of the residual symmetry under complex
        conjugation and $w\mapsto -w$. The lattice spacings are proportional
        to $n_{u}^{-1}, n_{w}^{-1}$. In both charts the precise domains
        of integration are moderately complicated, and so the exact number
        of points where the functions are evaluated to  approximate the integrals is not given by a
        simple formula. Writing $N_{1}(n_{x}, n_{p})$ for the number in
        the big chart and $N_{2}(n_{u}, n_{w})$ in the small chart, we have
       for example
        $$N_{1}(10,10)= 746,\ N_{1}(20,20)= 9644,\ N_{1}(30,30)=45481$$
        $$N_{2}(10,10)= 9737,\ N_{2}(20,20)=128149,\  N_{2}(24,20)=180435.
        $$

        \subsubsection{The volume form}
        
        To evaluate the volume form of the Fubini-Study metric defined
        by a given collection of co-efficients we work in the co-ordinates
        $x,y$. Let $V_{+}$ be the complex vector space of dimension $(k+1)(k+2)/2$
        with a basis labelled by the integer points in the \lq\lq big triangle''
        and $V_{-}$ the space of dimension $(k-2)(k-1)/2$ corresponding
        to the small triangle. Our data gives Hermitian metrics $\langle\
        ,\ \rangle_{\pm}$ on $V_{\pm}$
        and so a Hermitian metric $\langle\ ,\rangle$ on $V_{+}\oplus V_{-}$.
        Let $r(x,y)=r_{+}(x,y)\oplus r_{-}(x,y)$ be the vector valued function
        with the entries of $r_{+}(x,y)$ equal to $x^{p} y^{q}$ and the
        entries of $r_{-}(x,y)$ equal to $w x^{p} y^{q}$, where $w=\sqrt{1+x^{6}+y^{6}}$.
        (Of course $w$ is only defined up to a sign, but we choose a branch
        of the square root locally.) Write $r_{x}, r_{y}$ for the derivatives
        of $r$ with respect to $x,y$. Then away from the branch curve $w=0$
        the ratio $V=d\mu_{FS}/d\nu_{0}$
        of the Fubini-Study form and the fixed volume form $\theta\wedge
        \overline{\theta}$ is given by the formula
        \begin{equation}   V= \frac{1}{ \vert w\vert^{2} \Vert r\Vert^{6}}
        \rm{det} \left( \begin{array}{ccc} \langle r,r\rangle&\langle r,r_{x}\rangle&\langle
        r,r_{y}\rangle\\
        \langle r_{x},r\rangle&\langle r_{x},r_{x}\rangle & \langle r_{x},r_{y}\rangle 
        \\ \langle r_{y}, r\rangle& \langle r_{y}, r_{x}\rangle&\langle r_{y},
         r_{y}\rangle \end{array}\right) 
        \label{eq:determ}\end{equation}
        Of course this is the same whichever branch of the square root
        we take and the calculation is completely straightforward. The
        disadvantage is that we cannot use this formula near
         the branch curve because of the small denominator
          (although in practice the formula seems to be accurate
           close enough to the branch curve for most purposes).
             A better formula is as follows. Let $v_{\pm}$ be the vector
            -valued functions, taking values in $V_{\pm}$
             with entries $x^{p}y^{q}$. Thus $r_{+}=v_{+}$ and $r_{-}=w
             v_{-}$. Put 
             $$\delta_{x}=(v^{+}_{x}, w v^{-}_{x}), \delta_{y}=(v^{+}_{y},
             w v^{-}_{y})$$
             where subscripts denote differentiation. Set
             $$ \hdelx=\delta_{x}- \frac{\langle \delta_{x}, r\rangle}{\Vert
             r\Vert^{2}}\ \ , \ \ \hdely=\delta_{y}-\frac{\langle \delta_{y},
             r\rangle}{\Vert r \Vert^{2}}. $$
             Now write
             $$ f_{x}= 3 x^{5}, f_{y}=3 y^{5} . $$
             Define
             $$ Q_{x}= \frac{1}{\Vert r\Vert^{2}}\left( \langle \Vert v_{+}\Vert_{+}^{2}
             \langle v^{-}_{x}, v^{-}\rangle_{-}- \Vert v_{-}\Vert_{-}^{2}
             \langle v^{+}_{x}, v^{+}\rangle_{+}\right), $$
             and $Q_{y}$ symmetrically.
   Then $$V=\frac{1}{\Vert r\Vert^{4}}(V_{1}+V_{2}-V_{3}+ 2V_{4}-2V_{5}),
             $$
             where
             \begin{eqnarray*}   V_{1}&=& \vert w\vert^{2} \left( \Vert \hdelx\Vert^{2}
\Vert \hdely \Vert^{2}- \vert \langle \hdelx,\hdely\rangle\vert^{2}\right);
           \\
           V_{2}&=& \frac{\Vert v_{+}\Vert^{2}_{+}\Vert v_{-}\Vert^{2}_{-}}{\Vert
         r\Vert^{2}} \Vert f_{x} \hdely-f_{y} \hdelx \Vert^{2};\\
   V_{3}&=& \vert w\vert^{2} \vert f_{x} Q_{y}-f_{y} Q_{x} \vert^{2};\\
        V_{4}&=& {\rm Re}\left( \Vert \hdely \Vert^{2} \overline{Q_{x}}
      \overline{w}^{2} f_{x} + \Vert \hdelx \Vert^{2} \overline{Q_{y}}
           \overline{w}^{2} f_{y}\right); \\
     V_{5}&=& {\rm Re}\left( \langle \hdelx,\hdely\rangle ( \overline{f_{x}}
       Q_{y} w^{2} + f_{y} \overline{Q_{x}} \overline{w}^{2})\right).
           \end{eqnarray*}
            While it appears more complicated, this expression has advantages
            over Equation~\ref{eq:determ} even away from the branch curve.
            On the author's PC it takes about 1/5 of a second to evaluate
            this at a given point(when k=9). So to evaluate the volume form at
             30,000 
            points---which is the kind of number we use for our numerical
             integration---takes an hour or two.                     
  
\section{The Bergman kernel}

   Return to the general picture where $L\rightarrow X$ is an ample line
   bundle with a Hermitian metric and $d\nu$ is a fixed volume form on
   $X$. Then the space of sections $H^{0}(L^{k})$ has an $L^{2}$ Hermitian
   inner product. For a point $x\in X$ we have an evaluation map
   $$  e_{x}:H^{0}(L)\rightarrow L_{x}, $$
   which is represented by the inner product, so there is an element
   $\sigma_{x}\in H^{0}(L)\otimes \overline{L_{x}}$ such that
   $$     s(x) = \langle s, \sigma_{x}\rangle, $$
   for any section $s \in H^{0}(L)$. For $x,y\in X$ we define
   $$  K(x,y)= \vert \sigma_{x}(y)\vert^{2}, $$
   and the associated integral operator
   $$  Q_{K}(f)(x) = R \int_{X} K(x,y) f(y) d\nu_{y}, $$
   where $R={\rm Vol}(X,\nu)/{\rm dim}(H^{0}(L^{k})$. (This factor is included
   to make $Q_{K}$ independent of scalings of $\nu$, and $Q_{K}(1)=1$.)
   Let $s_{\alpha}$ be an orthonormal basis of $H^{0}(L^{k})$ with respect
   to the $L^{2}$ inner product. Then
   $$ \sigma_{x}(y)=\sum_{\alpha} s_{\alpha}(y) \otimes \overline{s_{\alpha}(x)}
   $$ and if we write
   $$f_{\alpha \beta}= (s_{\alpha}, s_{\beta}) $$
   (the pointwise inner product over $X$), we have
   $$ K(x,y)=\sum_{\alpha,\beta} f_{\alpha \beta}(x) f_{\beta\alpha}(y). $$
   Thus $Q_{K}$ is a finite-rank operator whose image lies in the finite-dimensional
   space $V\subset C^{\infty}(X)$ spanned by the $f_{\alpha \beta}$,{\it i.e.}
   the image $\iota({\cal H}(L^{k})$ in the notation of Section 2.2.1. 
   The
   restriction of $Q_{K}$ to $V$ gives an endomorphism of $V$ with 
   $$  Q_{K}(\sum a_{\gamma \delta} h_{\gamma \delta}) = R \sum a_{\gamma \delta}
   \langle f_{\gamma \delta}, f_{\beta\alpha}\rangle f_{\alpha \beta} $$
   where $\langle, \rangle$ denotes the $L^{2}$ inner product on functions.
   In other words we can define a endomorphism $Q:{\cal H}(L^{k})\rightarrow
   {\cal H}(L^{k})$ with matrix
      $$  Q_{\alpha \beta, \gamma \delta} = R \langle f_{\gamma \delta} f_{\alpha
   \beta}\rangle, $$
   and $\iota\circ Q= Q_{K}\circ \iota$.
   One interpretation of $Q$ is that it compares the two natural inner
   products on $V$. If we have a given metric on $H^{0}(L^{k})$ then we
   can identify ${\cal H}(L^{k})$ with the self-adjoint endomorphisms of
   $H^{0}(L^{k})$ and as such we have a standard  Hilbert-Schmidt norm on
    ${\cal H}(L^{k})$
   given by 
      $$ \Vert (a_{\alpha \beta})\Vert_{HS}^{2} = \sum \vert a_{\alpha
      \beta}\vert^{2}. $$
      On $V$ we have the restriction of the $L^{2}$ norm and these are
      related by
      $$   \Vert \iota(a) \Vert_{L^{2}}^{2}= R \langle a, Q(a) \rangle_{HS}.
      $$
      Notice that, if we regard ${\cal H}(L^{k})$ as the self-adjoint endomorphisms
      of ${\cal H}(L^{k})$ then $Q$ has been normalised so that $Q(1)=1$.
      
      To illustrate these ideas, take $L\rightarrow X$ to be the line bundle
      ${\cal O}(1)$ over $\bC\bP^{1}$ with the standard metric, and standard
      area form. As in Section (2.1.1) we take the usual $S^{1}$ action on $\bC\bP^{1}$
      and restrict attention to the invariant part of $V$, which corresponds
      to the diagonal matrices in ${\cal H}(L^{k}))$. With this restriction
      $Q$ is represented by a $(k+1)\times (k+1)$ matrix with entries
      \begin{equation}  Q_{ij} = \frac{k+1}{2k+1}  \ \frac{\left(\begin{array}{c}k\\i \end{array}\right)
      \left( \begin{array}{c} k\\j\end{array}\right)}{
      \left(\begin{array}{c}2k\\ i+j\end{array}\right)} 
      \ \ 0\,\leq i,j\leq k\label{eq:matrix} \end{equation}   
       (see the discussion in (4.3) below). 
       
       \subsection{The linearisation of the algorithm}
       
       One way in which the discussion above enters into our theory
       is in the analysis of the linearisation about a balanced metric.
        Fix an orthornormal basis $s_{\alpha}$ of
       $H^{0}(L^{k})$ for the balanced metric. Suppose $G_{\alpha \beta}=\delta_{\alpha
       \beta} + \epsilon_{\alpha \beta}$ is another metric. Then, to first
       order in $\epsilon$,
       $$  T_{\nu}(G)= R\int_{X} \frac{(s_{\alpha},s_{\beta})}{1- \sum \epsilon_{\alpha
       \beta} (s_{\alpha}, s_{\beta})} d\nu, $$
       so we have
       $$  T_{\nu}(G)= \delta_{\alpha \beta} + \tilde{\epsilon}_{\alpha
       \beta} + O(\epsilon^{2}), $$
       where 
       $$ \tilde{\epsilon}_{\alpha \beta}=R \int_{X} \sum_{\gamma \delta}
       (s_{\alpha}, s_{\beta})(s_{\gamma}, s_{\delta}) \epsilon_{\gamma
       \delta} d\nu, $$
    so $\tilde{\epsilon}=Q(\epsilon)$. In other words, the linearisation of the map $T$ at the balanced
       metric is given by $Q$. In particular, the largest eigenvalue of
       $Q$ on the trace-free matrices is the quantity $\sigma$ which determines
       the asymptotic rate of convergence of a sequence $T^{r}(G_{0})$,
       for almost all initial conditions $G_{0}$. Thus we can estimate
       this largest eigenvalue, in the examples discussed above, by analysing
       this convergence. For example, on the K3 surface $S$ we estimate,
       by analysing the sequences,
       that the eigenvalue is approximately $.22$ when $k=6$ and $.33$ when
       $k=9$.
       
      \subsection{Refined approximations and the heat kernel}
      
      Another way in which the operator $Q$ enters our theory is in the
      algorithm we have used for finding \lq\lq refined approximations''
      as discussed in (2.2.1) above. Recall that the linearisation of the map
      which assigns the volume form $\mu_{\omega}=\omega^{n}/n!$ to a Kahler metric
      $\omega$ is given by one half the Riemannian Laplacian, {\it i.e.}
      $$        \mu_{\omega+i \dbd  \phi}= \mu_{\omega}(1+ + \frac{1}{2}
      \Delta \phi) +O(\phi^{2}). $$
      Suppose given any metric $G$ on $H^{0}(L^{k})$, defining a metric
      $\omega$ on $X$, and let $s_{\alpha}$
      be an orthonormal basis of sections. Consider a small perturbation
      of $G$ to a metric with matrix $\delta_{\alpha \beta} +\epsilon_{\alpha
      \beta}$ in this basis. Then to first order in $\epsilon$ the induced
      Fubini-Study metric changes by $i \dbd  \phi$ where
      $$\phi= \sum \epsilon_{\alpha \beta} (s_{\alpha}, s_{\beta}).$$
      So, to first order, the change in the volume form is
      $$ \frac{1}{2} \Delta( \sum \epsilon_{\alpha \beta} (s_{\alpha},
      s_{\beta})).$$
      Now given a fixed volume form $\nu$ on $X$, write 
      $\mu_{\omega}=\eta
      \nu$, where we suppose $\eta$ is close to $1$. As in (2.2.1) define
      $$\eta_{\alpha\beta}= R \int_{X} (\eta-1) (s_{\alpha}, s_{\beta}) d\nu,
      $$
       and  consider the variation $$\epsilon_{\alpha \beta}= -\kappa
       \eta_{\alpha \beta}.$$
       Then the change in the volume form is, to first order,
       $$  \frac{\kappa}{2} \Delta Q_{K} (\eta). $$
       Thus the algorithm of (2.2.1) will replace an initial error term $\eta$
       by a new term which is approximately
       $$  W(\eta)= \eta-  \frac{\kappa}{2} \Delta Q_{K} (\eta), $$
       so we would like to argue that, for appropriate values of the parameter
       $\kappa$ and with respect to a suitable norm,
        the linear map $W$ is a {\it contraction}.
        
        To give evidence for this, we argue that the operator $Q_{K}$ should
        be related, asymptotically as $k\rightarrow \infty$, to the heat
        kernel on $X$. Consider the model case of a line bundle over $\bC^{n}$
        with curvature $-i\omega$, where $\omega$ is the standard Kahler
        form (corresponding to the Euclidean metric). Fix a trivialisation
        of the line bundle in which the connection form is $\frac{-i}{2}
        (\sum ( x_{a} dy_{a}- y_{a} dx_{a})$   where $z_{a}=x_{a}+ i y_{a}$
        are standard co-ordinates on $\bC^{n}$. Then, in this trivialisation,
         the section $\sigma_{0}$
        which represents evaluation at $0$ is
        $$ \sigma_{0}= \frac{1}{(2\pi)^{n}} e^{-\vert z\vert^{2}/4}, $$
        so our kernel is $$K(0,z)= \frac{1}{(2\pi)^{2n}} e^{-\vert z\vert^{2}/2}. $$
        The Euclidean heat kernel  is
        $$    H(0,z,t)= (\frac{1}{4\pi t})^{n}  e^{-\vert z\vert^{2}/4t}, $$
        so $K(0,z)=(2\pi)^{-n} H(0,z,1/2)$. Thus it is reasonable to expect that, on a general
        manifold $X$ the operator $Q_{K}$ will be approximately
         $e^{-\Delta/2}$ when $k$ is large; so the manifold has very large
         volume and the local geometry approaches the Euclidean case. (Notice
         that the factor of $2\pi^{n}$ is accounted for by the scaling
         built into the definition of $Q_{K}$, since when $k$ is large, by Riemann-Roch
         ${\rm dim} H^{0}(L^{k})$ is approximately $(2\pi)^{-n}$ times the
         volume, in the metric defined by the curvature form of $L^{k}$.)
         The author has not yet found any statement of exactly this kind
         in the literature but there are results very close to this in
         \cite{kn:DLM}, for example. In any case our present purpose is to give a
         plausible justification for the method rather than a rigorous
         proof. Of course, the Laplace operator considered above is that
         with respect to the \lq\lq large volume'' metric, with volume
         $O(k^{n})$.  
         With this discussion in place, we argue that near to the Calabi-Yau
         metric, the operator $W$ is approximately
         $$ \tilde{W}= 1- \kappa \frac{\Delta}{2} e^{-\Delta/2}. $$
         Now $\tilde{W}$ is easy to analyse in terms of the spectrum 
         of the Laplacian. On an eigenspace belonging to eigenvalue $\mu$
         $\tilde{W}$ acts as $(1-\kappa \frac{\mu}{2} e^{-\mu/2})$. Since
         the function $xe^{-x}$ has maximum value $e^{-1}$ for positive
         $x$, the operator $\tilde{W}$ is a contraction provided that
         $ 0< \kappa<2e $. This is consistent with the values of the parameter
         $\kappa$ found to be effective empirically. These ideas also
          explain why
         the \lq\lq refining algorithm'' takes a long time to get very
         close to the refined approximation, since the contraction factor
         for large eigenvalues $\mu$ is extremely close to $1$.
         
         As a byproduct of these ideas, we can hope to get information
         about the spectrum of the Laplacian of the Calabi-Yau metric from
         our theory. Let $\Delta_{0}$ be the Laplacian of the metric scaled
         to have total volume $(2\pi)^{n}$ and write
         $$  k'= ({\rm dim}\ H^{0}(L^{k}))^{1/n}. $$
          Then  we expect that the spectrum of
         $Q$ approximates that of 
         $   e^{-\Delta_{0}/2k'} $. 
         Thus if $\lambda$ is the first eigenvalue of $\Delta_{0}$ we expect
         that the convergence parameter $\sigma_{k}$ associated to our
         algorithm is approximately            
         $ e^{-\lambda/2k'} $. 
          If, as in this paper, we work with $\Gamma$-invariant metrics
          then we should take $\lambda$ to
         be the first eigenvalue on $\Gamma$-invariant functions. Our estimates
         $\sigma_{6}=.22,\sigma_{9}= .33$ are reasonably consistent with
         this  since
         $$ -2 \log(.22) (38)^{1/2}= 18.7$$
         and $$-2\log(.33) (83)^{1/2}= 20.2.$$
         So we expect that $\lambda$ is about 20. 
          (We can also hope to get approximations to the eigenfunctions
         of the Laplacian from the eigenvectors of Q.)   
         
         \subsection{Algebraic approximation to the heat kernel.}
         
         We have now explained the importance of the finite-dimensional
         linear operator $Q$ in our theory, and its (probable) relation to the
         Laplace operator on the manifold. Recall that the matrix entries
         of $Q$, in terms of an orthonormal basis $s_{\alpha}$ of $H^{0}(L^{k})$,
         are
         $$  Q_{\alpha \beta, \gamma \delta} = R \int_{X}  (s_{\alpha},
         s_{\beta})(s_{\gamma}, s_{\delta}) d\nu. $$
         On the face of it, this requires us to evaluate the large number
         ${\rm dim} H^{0}(L^{k})^{4}$ of integrals over $X$ to find the matrix. However
         we can write
         $$  Q_{\alpha \beta, \gamma \delta}= \int_{X} (s_{\alpha}
         s_{\delta}, s_{\beta}s_{\gamma}) d\nu, $$
         where the products are sections of $L^{2k}$ and $(\ ,\ )$ denotes
         the fibre metric on $L^{2k}$. Let $\tau_{i}$ be a basis of $H^{0}(L^{2k})$.
         If we know the integrals
         $$ I_{ij}= \int_{X} (\tau_{i}, \tau_{j}) d\nu, $$
         then we can compute the matrix entries in terms of purely {\it
         algebro-geometric} data: the product map
         \begin{equation} H^{0}(L^{k}) \otimes H^{0}(L^{k}) \rightarrow H^{0}(L^{2k}).
         \label{eq:square}\end{equation} Explicitly, if
         $$ s_{\alpha} s_{\beta} = \sum P_{\alpha \beta i} \tau_{i}$$ then
         $$ L_{\alpha \beta, \gamma \delta}= \sum_{ij} P_{\alpha \delta i} \overline{P_{\beta
         \gamma j}} I_{ij} . $$
         This means that we only need to evaluate approximately $ 2^{n}
         (N+1)^{2}$ integrals to find the matrix. Moreover these integrals
         are precisely the integrals which define the map $T$ for the line
         bundle $L^{2k}$. In geometric terms, for any vector space $V$
         we have the Veronese embedding
         $$ \bP(V)\rightarrow \bP(s^{2} (V)). $$
         A hermitian metric on $V$ defines a standard induced metric on
         $s^{2} V$ and, up to a scale factor, the Veronese embedding is
         isometric with respect to the Fubini-Study metrics. Thus for $X$
         in $\bP(V)$ we get the same induced metric by embedding
         in $\bP(s^{2}(V))$. Starting with the canonical embedding in
          $V=H^{0}(L^{k})^{*}$ we get the canonical embedding in $\bP
          (H^{0}(L^{2k}))^{*}$,
          which is contained as a linear subspace in $\bP(s^{2} V)$.
          Thus, starting with a metric $G$ on $H^{0}(L^{k})$ we take the standard
          induced metric $G'$ on $H^{0}(L^{2k})$, regarded as a quotient of
          the symmetric product. Then the calculation of $T(G')$ is equivalent
          to the calculation of the matrix entries. 
          
          Now suppose it happens that $G$ is the balanced metric for $L^{k}$
          and $G'$ is also the balanced metric for $L^{2k}$. This will
          only be the case in rather special circumstances, but for example
          it holds when $X=\bC \bP^{1}$ with $\nu$ equal to the standard
          area form. In this case we have $T(G')=G'$ and we can find the
          matrix entries purely algebraically, in terms of the product
          map Equation~\ref{eq:square} and the original hermitian metric $G$. For example for
          the line bundle ${\cal O}(k)$ over $\bC\bP^{1}$, restricting
          to the $S^{1}$-invariant metrics, we get the matrix entries $Q_{ij}$
          in Equation~\ref{eq:matrix} above. But in any case we can {\it define} another
          endomorphism $\tilde{Q}$ of ${\cal H}(L^{k})$ by this procedure. That
          is, we take the matrix enties
          $$ \tilde{Q}_{\alpha \beta} = \sum_{ij} P_{\alpha \delta i}
          \overline{P}_{\beta \gamma j}\tilde{I}_{ij},
          $$ where $\tilde{I}_{ij}$ are the inner products in $H^{0}(L^{2k})$
          given by the induced hermitian metric, regarded as a quotient
          of the symmetric square. It is reasonable to expect that, when
          $k$ is sufficiently large, the metric $G'$ is close to the balanced
          metric and hence that $\tilde{Q}$ is a good approximation to
          $Q$. To sum up, starting with a hermitian metric $G$ on $H^{0}(L^{k})$
          we have a purely algebraic procedure for defining a self-adjoint
          endomorphism $\tilde{Q}$ on ${\cal H}(L^{k})$, and when $G$ is the balanced
          metric (or close to the balanced metric) we can expect that
          $\tilde{Q}$ is an approximation to the heat kernel
          $e^{-\Delta_{0}/2k'}$.
          
          To illustrate these ideas consider first the case of $S^{1}$-invariant
          metrics on $S^{2}$. The $SU(2)$ invariance of the problem
          implies that the eigenspaces of $Q_{K}$ correspond to spherical harmonics.
          Let $z\in [-1,1]$ be the standard height co-ordinate on the
           sphere and $p\in S^{2}$ be the pole where $z=1$. The kernel
           function $K(p,\ )$ associated to ${\cal O}(k)$ is proportional
           to $(1+z)^{k}$ and so the eigenvalue $\chi_{m,k}$ of $Q$ associated
           to the spherical harmonics of degreee $m$ is
           $$  \chi_{m,k}= \frac{k+1}{2^{k+1}}
            \int_{-1}^{1} (1+z)^{k} P_{m}(z) dz, $$
            where $P_{m}$ is the Legendre polynomial. It is an exercise
            in Legendre polynomials to show that
            $$ \chi_{m,k}= \frac{(k'-1) \dots (k'-m)}{(k'+1) \dots (k'+m)},
            $$
            where we write $k'= k+1= {\rm dim } H^{0}({\cal O}(k))$, 
            as above. 
              Now set $\lambda_{m,k}= - 2 k' \log(
          \chi_{m,k})$, so
          $$   \lambda_{m,k}= -2k' \sum_{r=1}^{m} \left(\log(1-\frac{r}{k'}) -
          \log(1+\frac{r}{k'})\right) . $$
          From the Taylor expansion of the logarithm we  see that
           $$  \lambda_{m,k} = 2m (m+1) + O(k^{-2}), $$
           and the limits  $2m(m+1)$ are the eigenvalues of the Laplacian
           on the sphere of area $2\pi$. For example, we have
           $$ \lambda_{1,4}= 4+.055, \lambda_{2,10}= 12+.15 , \lambda_{4,30}=
           40+.14 . $$
          
          This also ties in with observed value of $\sigma$ in Section
          (2.2), since $\chi_{2,6}=5/12=.41666\dots $.
          
          Finally we consider the balanced metric on the $K3$ surface $S$,
          with $k=6$. We restrict to the $\Gamma$-invariant part of ${\cal
          H}(L^{k})$. We work using  a orthonormal basis of $H^{0}(L^{k})$ given by rescaling
          the standard monomials, apart from the triple $1,x^{6}, y^{6}$.
          Here we choose scalars $A,B$ such that $A+Bx^{6}+By^{6}$ and
          the two similar terms given by permutations are orthonormal.
          Then we find that the endomorphism $\tilde{Q}$ on the $11$-dimensional
           $\Gamma$-invariant
          part of ${\cal H}(L^{k})$ has matrix

  $$ 10^{-2} \times \left(\begin{array}{ccccccccccc}
    .61&-.33&-.81&-.38&-.84&-.91&-.10&.96&.36&1.19&1.45\\
    \ &3.68&8.38&5.43&4.09&5.71&5.15&2.30&3.13&5.44&2.39\\
    \ &\ &19.6&13.0&10.6&14.8&13.5&6.46&7.50&13.9&6.56\\
    \ &\ &\ &8.96&7.63&10.7&10.2&5.63&5.27&10.5&5.72\\
    \ &\ &\ &\ &9.21&12.7&11.5&6.32&4.34&10.8&6.04\\
    \ &\ &\ &\ &\ &17.7&16.3&9.62&6.23&15.7&9.35\\
    \ &\ &\ &\ &\ &\ &16.9&12.5&6.38&16.9&12.4\\
    \ &\ &* &\ &\ &\ &\ &14.4&4.01&12.9&14.4\\
    \ &\ &\ &\ &\ &\ &\ &\ &3.47&7.29&4.35\\
    \ &\ &\ &\ &\ &\ &\ &\ &\ &19.0&13.7\\
    \ &\ &\ &\ &\ &\ &\ &\ &\ &\ &14.7 
    \end{array}\right) $$ 
      The matrix is symmetric so we omit the entries below the diagonal.
      Here the first basis element corresponds to the off-diagonal term,
   the next seven to the entries in the big triangle and the last three
   to the small triangle.       
          
          We find the first six eigenvalues (ordered by absolute value)
           of this matrix numerically. They are
          $$1.002,.1956,.05857,.02395,-.002669,.002388 . $$
          The first eigenvalue, 1.002, is a substitute for the exact eigenvalue
          1 of the matrix $Q$, and the close agreement is encouraging. The
          fact that the fifth eigenvalue is negative, whereas $Q$ is a
          positive operator, shows that we cannot take the approximation
          this far
          down the spectrum. For each positive eigenvalue $\chi$ we compute
          $\lambda= - 2\sqrt{38} \log \chi$ and we for the second, third 
          and
          fourth eigenvalues we obtain the $\lambda$-values
          $$  20.12, 34.98, 46.00 $$
          respectively.
          The  eigenvalue $\chi=.1956$ is  in reasonable agreement with 
          our previous numerical estimate .22 for the first eigenvalue of $Q$ and the
          corresponding estimate $20.12$ for the lowest eigenvalue of the
          Laplacian is very close to our previous estimate $20.2$ from
          the observed value of $\sigma $ when  
          case $k=9$. 
          It is perhaps reasonable to predict, based on this discussion,
           that the next eigenvalue of
          the Laplacian  $\Delta_{0}$ (on invariant functions) is  about $35$.
           It would be interesting to test this by repeating
          the work for $k=9$, but the author has not yet had time to do
          so.

\section{Appendix}

Here we discuss the fact stated in the intoduction, that given any Kahler
metric $\omega$ in the class $c_{1}(L)$ there is a sequence of \lq\lq algebraic
'' metrics $\omega_{k}$ arising from Hermitian forms on $H^{0}(L^{k})$,
with $\omega_{k}-\omega=o(k^{\nu})$ for all $\nu$. 

The proof uses the Tian-Zelditch-Lu expansion and is similar to the argument
in \cite{kn:D1}.  We start with Tian's 
approximation
 which, in the notation of Section 2, 
  is to take $ FS\circ {\rm Hilb}(\omega)$. (Here we are regarding
  $k$ as a parameter which is supressed in the notation.) Then
  $$ FS\circ {\rm Hilb} (\omega) = \omega+ k^{-1} i \dbd  \log(\rho_{\omega})$$
  where $\rho_{\omega}$ is the density of states function $\sum \vert s_{\alpha}
  \vert^{2}$ for an orthonormal basis $s_{\alpha}$. We know that
  $\rho_{\omega}$ has an asymptotic expansion
  $$\rho_{\omega}= 1+ k^{-1} a_{1}(\omega)+ k^{-2} a_{2}(\omega) + \dots,
  $$ for certain local invariants $a_{i}$ of the Kahler metric $\omega$.
  Thus $\tilde{\omega}-\omega=O(k^{-2})$ and the order $k^{-2}$ term is
  $i \dbd  a_{1}(\omega)$. Now let $\omega^{*}= \omega-k^{-2}i \dbd
  a_{1}(\omega)$ and consider the metric $FS\circ Hilb(\omega^{*})$. Applying
  the expansion, with smooth dependence on parameters, we see that 
  $$ FS\circ {\rm Hilb} (\omega^{*})= \omega + O(k^{-3}). $$
  We can repeat this process to kill of successively as many terms as we
  please in the asymptotic expansion. The correction terms will become
  progressively more complicated, involving contributions from the derivatives
  of $a_{i}(\omega)$ with respect to $\omega$, just as in \cite{kn:D1}. In this way
  we obtain, for any given $\nu$, a seqence of approximations $\omega_{k}$
  with  $\omega_{k}=\omega + o(k^{\nu})$. A standard \lq\lq diagonal'
  argument gives a single sequence with difference $o(k^{\nu})$ for any
  $\nu$.

  There is an elementary argument to prove a somewhat weaker result.
   The construction of a Fubini-Study metric $\omega_{H}$ from a Hermitian
  form $H$ can be extended to allow {\it indefinite} forms $H$, so long
  as $H$ is positive on the vectors in $H^{0}(L^{k})^{*}$ corresponding to
  points of $X$. It is easy to prove that any Kahler metric can be rapidly
  approximated by algebraic metrics in this larger class.  It is convenient to assume
   that $L$ is a
  {\it very ample} line bundle over $X$, so the sections of $L$ give an
  embedding of $X$ in $\bC\bP^{N}$. (The argument can be extended to avoid
  this assumption.) For $k\geq 1$ consider the standard Veronese embedding
  $$  \bC\bP^{N}\rightarrow \bC\bP^{N_{k}}. $$
  Let $Z_{A}$ denote standard homogeneous co-ordinates on $\bC\bP^{N_{k}}$
  and let $V_{k}\subset C^{\infty}(\bC\bP^{N})$ be the vector space of
  complex-valued functions on $\bC\bP^{N}$ given by linear combinations of
   $$ \frac{Z_{A}\overline{Z_{B}}}{\vert z\vert^{2k}} $$
   Of course the $Z_{A}$ are just the monomials of degree $k$ in the homogeneous
   co-ordinates $z_{\alpha}$ on $\bC\bP^{N}$. 
   \begin{lem}
   The space $V_{k}$ is the direct sum of the first $k$ eigenspaces of
   the Laplace operator $\Delta_{\bC\bP^{N}}$ for the standard Fubini-Study
   metric on $\bC\bP^{N}$.
   \end{lem}
   To prove this we take as known the analogous and well-known fact for the Laplacians
   on spheres. The sum of the first $k$ eigenspaces for the Laplacian on
    $S^{m-1}\subset \bR^{m}$ is exactly the space of functions given by
    restrictions of polynomials of degree $k$ on $\bR^{m}$.  
   Now consider the Hopf fibration $S^{2N+1}\rightarrow \bC\bP^{N}$. This
   is a Riemannian submersion so eigenfunctions of the Laplacian on $\bC\bP^{N}$
   lift to $S^{1}$-invariant eigenfunctions on the sphere. So the sum of
   the first $k$ eigenspaces on $\bC\bP^{N}$ can be identified with the
   polyomials in $z_{\alpha},\overline{z_{\alpha}}$ which are $S^{1}$-invariant.
   But it is clear that these are just polynomials in the products
    $z_{\alpha}\overline{z_{\beta}}$. Separating out the holomorphic and
    antiholomorphic terms, we see that these are exactly the linear combinations of
    the products $Z_{A}\overline{Z_{B}}$, as required.   
    
    Now it is a standard fact that if $f$ is a smooth function on a compact
    Riemannian manifold and $f_{k}$ is the $L^{2}$ projection of $f$ to the sum
    of the first $k$ eigenspaces of the Laplacian then $f_{k}-f= o(k^{\nu})$
    for any $\nu$. Let $\omega_{0}$ be the metric on $X$ given by the
    restriction of the standard Fubini-Study metric on $\bC\bP^{N}$ so
    $$ \omega=\omega_{0}+i\dbd \phi, $$
    for some  smooth function $\phi$ on $X$. Extend $\phi$ arbitrarily
    to a smooth function on $\bC\bP^{N}$ and take $f=e^{\phi}$. Then
    $f$ is a positive real valued function on $\bC\bP^{N}$ and there is
    no loss in generality in supposing that the projections $f_{k}$ are also positive
    on $\bC\bP^{N}$. By the Lemma the function $f_{k}$ is a sum
    $$ \frac{ \sum h_{A B} Z_{A} \overline{Z_{B}}}{\vert z \vert^{2k}}.
    $$ and the Fubini-Study metric $\omega_{k}$ associated to the form with matrix
    $\delta_{AB}+h_{AB}$ is $\omega +i \dbd (\log f_{k}- \log f) $. So
    $\omega_{k}-\omega$ is $o(k^{\nu}) $ for all $\nu$.


\begin{thebibliography}{99}
     \bibitem{kn:A} D. Acheson {\em From calculus to chaos} Oxford U.P.
     1997
     \bibitem{kn:Ar} C. Arezzo, A. Ghigi and A. Loi {\em Stable bundles and the first eigenvalue of the
     Laplacian} Preprint (2005) 
  \bibitem{kn:BLY} J-P. Bourguignon, P.Li and S-T. Yau {\em Upper bounds
  for the first eigenvalue of algebraic submanifolds} Comment. Math. Helvetici
  69 (1994) 199-207
  \bibitem{kn:C} E. Calabi {\em Extremal Kahler metrics} In:Seminar on
  Differential Geometry S-T. Yau ed., Princeton U.P. 1983
  \bibitem{kn:DLM} X Dai, K. Liu and X. Ma {\em On the asymptotic expansion
  of the Bergman kernel} Preprint
  \bibitem{kn:D1} S. Donaldson {\em Scalar curvature and projective embeddings,
  I} Jour. Differential Geometry 58 (2001) 479-522 
  \bibitem{kn:D2} S. Donaldson {\em Scalar curvature and projective embeddings,
  II} Quarterly Jour. Math.  56 (2005) 345-56
\bibitem{kn:HW} M. Headrick and T. Wiseman {\em Numerical Ricci-flat metrics
on K3} hep-th/0506129
\bibitem{kn:Luo} H. Luo {\em Geometric criterion for Mumford-Gieseker stability
of polarised manifold} Jor. Differential Geometry 49 (1998) 577-99
\bibitem{kn:Mab} T. Mabuchi {\em An energy-theoretic approach to the Hitchin-Kobayashi
correspondence for manifolds} Invent. Math. 159 (2004) 225-243
\bibitem{kn:MZ} J. Millson and B. Zombro {\em A Kahler structure on the
moduli space of isometric maps of a circle into Euclidean space}
Invent. Math. 123 (1996) 35-59
\bibitem{kn:Mum} D. Mumford {\em Geometric Invariant Theory} Springer 1982
\bibitem{kn:New} P. Newstead {\em Introduction to moduli problems and orbit
spaces} Tata Institute Lectures 51, Springer 1978 
\bibitem{kn:S} Y. Sano {\em Numerical algorithm for finding balanced metrics}
Tokyo Institute of Technology Preprint, 2004 
\bibitem{kn:T} G. Tian {\em On a set of polarised Kahler metrics on algebraic
manifolds} J. Differential Geometry 32 (1990) 99-130
\bibitem{kn:Zh} S. Zhang {\em  Heights and reductions of semistable varieties}
Compositio Math. 104 (1996) 77-105

\end{thebibliography}
\end{document}